\documentclass[11pt,twoside]{article}
\newcommand{\ifims}[2]{#1} 
\newcommand{\ifAMS}[2]{#1}   
\newcommand{\ifau}[4]{#1}  
\newcommand{\ifbook}[2]{#1}   

\usepackage{ifthen}
\usepackage[ampersand]{easylist}
\usepackage{amsmath,amssymb,amsthm}
\usepackage[bb=stix]{mathalpha}
\usepackage{mathtools}
\usepackage{natbib}
\usepackage{epsfig,graphicx}
\usepackage{comment}
\usepackage{color}
\usepackage{srcltx}
\usepackage[mathscr]{eucal}
\usepackage[math]{easyeqn}
\usepackage{etoolbox}
\usepackage{hyperref}
\hypersetup{
            colorlinks,
            linkcolor=hookersgreen,
            linktoc=blue,
            citecolor=ultramarine,
            urlcolor=black,
            filecolor=black
            }
\usepackage{nameref}
\usepackage{multirow}

\usepackage{accents}
\usepackage{mathrsfs}
\usepackage{bbm}
\usepackage{algorithm}
\usepackage{algpseudocode}

\ifims{
\textheight=22cm
\textwidth=14.8cm
\topmargin=0pt
\oddsidemargin=1.0cm
\evensidemargin=1.0cm
\linespread{1.3}

}{ 
}

\numberwithin{equation}{section}
\numberwithin{figure}{section}
\newcounter{example}[section]
\numberwithin{example}{section}
\newcounter{remark}[section]
\numberwithin{remark}{section}
\newtheorem{theorem}{Theorem}[section]
\newtheorem{proposition}[theorem]{Proposition}
\newtheorem{lemma}[theorem]{Lemma}
\newtheorem{corollary}[theorem]{Corollary}

\newtheorem{exmp}[example]{Example}
\newtheorem{rmrk}[remark]{Remark}
\newenvironment{example}{\begin{exmp}\rm}{\end{exmp}}
\newenvironment{remark}{\begin{rmrk}\rm}{\end{rmrk}}

\ifbook{
    
    \newcommand{\Section}[1]{\subsection{#1}}
    \newcommand{\Subsection}[1]{\subsubsection{#1}}

  }
  {
    
    \newcommand{\Section}[1]{\section{#1}}
    \newcommand{\Subsection}[1]{\subsection{#1}}

  }

    \renewcommand{\Section}[1]{\section{#1}}
    \renewcommand{\Subsection}[1]{\subsection{#1}}

\renewcommand{\(}{$\,}
\renewcommand{\)}{\,$}

\def\nquad{\hspace{-1cm}}
\def\eqdef{\stackrel{\operatorname{def}}{=}}
\def\eqd{\stackrel{\operatorname{d}}{=}}

\DeclareMathAlphabet{\mathbbmsl}{U}{bbm}{bx}{sl}

\DeclareMathSymbol{\Alpha}{\mathalpha}{operators}{"41}
\DeclareMathSymbol{\Beta}{\mathalpha}{operators}{"42}
\DeclareMathSymbol{\Epsilon}{\mathalpha}{operators}{"45}
\DeclareMathSymbol{\Zeta}{\mathalpha}{operators}{"5A}
\DeclareMathSymbol{\Eta}{\mathalpha}{operators}{"48}
\DeclareMathSymbol{\Iota}{\mathalpha}{operators}{"49}
\DeclareMathSymbol{\Kappa}{\mathalpha}{operators}{"4B}
\DeclareMathSymbol{\Mu}{\mathalpha}{operators}{"4D}
\DeclareMathSymbol{\Nu}{\mathalpha}{operators}{"4E}
\DeclareMathSymbol{\Omicron}{\mathalpha}{operators}{"4F}
\DeclareMathSymbol{\Rho}{\mathalpha}{operators}{"50}
\DeclareMathSymbol{\Tau}{\mathalpha}{operators}{"54}
\DeclareMathSymbol{\Chi}{\mathalpha}{operators}{"58}
\DeclareMathSymbol{\omicron}{\mathord}{letters}{"6F}

\newcommand{\cc}[1]{\mathscr{#1}}
\newcommand{\bb}[1]{\boldsymbol{#1}}

\renewcommand{\bar}[1]{%
  \hbox{%
    \vbox{%
      \hrule height 0.5pt 
      \kern0.35ex
      \hbox{%
        \kern-0.05em
        \ensuremath{#1}%
        \kern-0.05em
      }%
    }%
  }%
} 

\makeatletter
\newcommand*\rel@kern[1]{\kern#1\dimexpr\macc@kerna}
\newcommand*\widebar[1]{%
  \begingroup
  \def\mathaccent##1##2{%
    \rel@kern{0.8}%
    \overline{\rel@kern{-0.8}\macc@nucleus\rel@kern{0.2}}%
    \rel@kern{-0.2}%
  }%
  \macc@depth\@ne
  \let\math@bgroup\@empty \let\math@egroup\macc@set@skewchar
  \mathsurround\z@ \frozen@everymath{\mathgroup\macc@group\relax}%
  \macc@set@skewchar\relax
  \let\mathaccentV\macc@nested@a
  \macc@nested@a\relax111{#1}%
  \endgroup
}
\makeatother

\renewcommand{\hat}[1]{\widehat{#1}}
\renewcommand{\tilde}[1]{\widetilde{#1}}

\makeatletter
\def\mathcenterto#1#2{\mathclap{\phantom{#1}\mathclap{#2}}\phantom{#1}}
\let\old@widetilde\widetilde
\def\widetildeto#1#2{\mathcenterto{#2}{\old@widetilde{\mathcenterto{#1}{#2\,}}}}
\let\old@widehat\widehat
\def\widehatto#1#2{\mathcenterto{#2}{\old@widehat{\mathcenterto{#1}{#2\,}}}}
\makeatother

\newcommand{\thankstitle}[1]{\ifthenelse{\equal{#1}{}}{}{\thanks{#1}}}
\newcommand{\thanksau}[1]{\ifthenelse{\equal{#1}{}}{}{\thanks{#1}}}

\newcommand{\aua}[6]
{\def\authora{#1}
\def\runauthora{#2}
\def\addressa{#3}
\def\emaila{#4}
\def\affiliationa{#5}
\def\thanksa{#6}}

\def\theauthors{
\ifau{ 
  \author{
    \authora
    \thanksau{\thanksa}
    \\[5.pt]
    \addressa \\
    \texttt{ \emaila}
  }
}
{  
  \author{
    \authora
    \thanksau{\thanksa}
    \\[5.pt]
    \addressa \\
    \texttt{ \emaila}
    \and
    \authorb
    \thanksau{\thanksb}
    \\[5.pt]
    \addressb \\
    \texttt{ \emailb}
  }
}
{   
  \author{
    \authora
    \thanksau{\thanksa}
    \\[5.pt]
    \addressa \\
    \texttt{ \emaila}
    \and
    \authorb
    \thanksau{\thanksb}
    \\[5.pt]
    \addressb \\
    \texttt{ \emailb}
    \and
    \authorc
    \thanksau{\thanksc}
    \\[5.pt]
    \addressc \\
    \texttt{ \emailc}
  }
} {   
  \author{
    \authora
    \thanksau{\thanksa}
    \\[5.pt]
    \addressa \\
    \texttt{ \emaila}
    \and
    \authorb
    \thanksau{\thanksb}
    \\[5.pt]
    \addressb \\
    \texttt{ \emailb}
    \and
    \authorc
    \thanksau{\thanksc}
    \\[5.pt]
    \addressc \\
    \texttt{ \emailc}
    \and
    \authord
    \thanksau{\thanksd}
    \\[5.pt]
    \addressd \\
    \texttt{ \emaild}
  }
}
}

\renewcommand{\Gamma}{\varGamma}
\renewcommand{\Pi}{\varPi}
\renewcommand{\Sigma}{\varSigma}
\renewcommand{\Delta}{\varDelta}
\renewcommand{\Lambda}{\varLambda}
\renewcommand{\Psi}{\varPsi}
\renewcommand{\Phi}{\varPhi}
\renewcommand{\Theta}{\varTheta}
\renewcommand{\Omega}{\varOmega}
\renewcommand{\Xi}{\varXi}
\renewcommand{\Upsilon}{\varUpsilon}

\def\av{\bb{a}}

\def\uv{\bb{u}}

\def\xv{\bb{x}}

\def\Av{\bb{A}}

\def\Xv{\bb{X}}
\def\Yv{\bb{Y}}

\def\epsv{\bb{\varepsilon}}

\def\xiv{\bb{\xi}}

\def\Psiv{\bb{\Psi}}

\def\sumi{\sum_{i=1}^{n}}

\usepackage{color}

\definecolor{blue(pigment)}{rgb}{0.2, 0.2, 0.6}
\definecolor{ultramarine}{rgb}{0.07, 0.04, 0.56}
\definecolor{darkspringgreen}{rgb}{0.09, 0.45, 0.27}
\definecolor{hookersgreen}{rgb}{0.0, 0.44, 0.0}
\definecolor{plum(traditional)}{rgb}{0.56, 0.27, 0.52}
\definecolor{purple(html/css)}{rgb}{0.5, 0.0, 0.5}
\definecolor{magenta(dye)}{rgb}{0.79, 0.08, 0.48}

\def\errS{\mathcal{E}}
\def\errSt{\tilde{\errS}}

\def\constg{\kappa}
\def\nuH{\nu}

\def\Pmuvp{\Phi}
\def\zzH{\zz}

\def\muHa{\muH_{\alpH}}
\def\muHc{\muH_{c}}

\def\KH{\KS}

\def\BHT{\mathbbmss{B}}

\def\alpH{\alpha}

\def\Gauss{\mathcal{G}}
\def\Gaussb{\bar{\Gauss}}

\def\Sigmah{\hat{\Sigma}}

\def\rhoH{\rho}

\def\supH{\lambda}
\def\rexH{\epsilon}

\def\Egs{\E_{\gaussv}}
\def\Pgs{\P_{\gaussv}}

\def\Matr{\mathfrak{M}}

\def\weights{\weight^{*}}

\def\vH{\vA}

\def\Eta{\mathcal{H}}

\def\nbin{N}

\def\HVB{\mathbbmsl{V}} 

\def\dltwu{\tau}

\def\II{\mathcal{I}}

\def\R{\mathbbmsl{R}}
\def\E{\mathbbmsl{E}}

\def\P{\mathbbmsl{P}}

\def\kappa{\varkappa}

\def\Bernoulli{\operatorname{Bernoulli}}

\def\diag{\operatorname{diag}}

\def\Fr{\operatorname{Fr}}

\def\ND{\mathcal{N}}

\def\oper{\operatorname{op}}

\def\Var{\operatorname{Var}}

\def\T{\top}
\def\tr{\operatorname{tr}}

\def\nsize{{n}}

\def\sumi{\sum_{i=1}^{\nsize}}

\def\ex{\mathrm{e}}

\def\Id{\mathbbmsl{I}}
\def\Ind{\operatorname{1}\hspace{-4.3pt}\operatorname{I}}

\def\alp{\alpha}

\def\BBH{B}

\def\cdens{\phi}

\def\CONST{\mathtt{C} \hspace{0.1em}}
\def\CONSTi{\mathtt{C}}

\def\dimA{\mathtt{p}}

\def\dimH{\dimA}

\def\dimp{p}

\def\dimq{q}

\def\Eta{\cc{H}}

\def\eps{\epsilon}			
\def\eps{\varepsilon}

\def\gaussv{\bb{\gauss}}
\def\gauss{\gamma}

\def\gm{\mathtt{g}}

\def\gmb{\gm}

\def\gmn{\gm}
\def\imi{\mathtt{i}}

\def\KS{A}

\def\Kappa{\cc{K}}

\def\muH{\mu}

\def\norms{\circ} 

\def\QP{Q}

\def\thetas{\theta^{*}}

\def\Tau{T}

\def\ups{\upsilon}

\def\ups{\upsilon}

\def\upss{\ups^{*}}

\def\vA{\mathtt{v}}

\def\VP{V}

\def\weight{w}

\def\xx{\mathtt{x}}
\def\xxc{\xx_{c}}

\def\xxs{\xx_{\norms}}

\def\xxe{\xx_{\ex}}

\def\yy{\mathtt{y}}

\def\zq{z}
\def\zqc{\zq_{c}}

\def\zqs{\bar{\zq}}

\def\zqe{\mathsf{z}}

\def\zz{\mathfrak{z}}

\def\thetitle{Deviation bounds for the norm of a random vector under exponential moment conditions with applications}
\def\theruntitle{Deviation bounds under exponential moment conditions}

\def\theabstract{
Hanson-Wright inequality provides a powerful tool for bounding the norm \( \| \xiv \| \) of a centered stochastic vector \( \xiv \)
with sub-gaussian behavior. 
This paper extends the bounds to the case when \( \xiv \) only has bounded exponential moments of the form
\( \log E \exp \langle V^{-1} \xiv,u \rangle \leq \| u \|^{2}/2 \), where \( V^{2} \geq \mathrm{Var}(\xiv) \) and \( \| u \| \leq g \) for some fixed \( g \).
For a linear mapping \( Q \), we present an upper quantile function \( z_{c}(B,x) \) ensuring 
\( P(\| Q \xiv \| > z_{c}(B,x)) \leq 3 e^{-x} \) with \( B = Q \, V^{2} Q^{T} \).
The obtained results exhibit a phase transition effect: with a value \( x_{c} \) depending on \( g \) and \( B \),
for \( x \leq x_{c} \), 
the function \( z_{c}(B,x) \) replicates the case of a Gaussian vector \( \xiv \), that is,
\( z_{c}^{2}(B,x) = \tr(B) + 2 \sqrt{x \tr(B^{2})} + 2 x \| B \| \).
For \( x > x_{c} \), the function \( z_{c}(B,x) \) grows linearly in \( x \).
The results are specified to the case of Bernoulli vector sums and to covariance estimation in Frobenius norm.
}

\def\kwdp{62E15}
\def\kwds{62E10}

\def\thekeywords{upper quantiles, phase transition, vector Bernoulli sums, Frobenius loss}

\def\thankstitle{}

\aua
{Vladimir Spokoiny}
{Spokoiny, V. }
{
    Weierstrass Institute and HU Berlin,  
    \\
    Mohrenstr. 39, 10117 Berlin, Germany
    }
{spokoiny@wias-berlin.de}
{Weierstrass-Institute and Humboldt University Berlin}
{
  Financial support by the German Research Foundation (DFG) through the Collaborative Research Center 1294 ``Data assimilation'' is gratefully acknowledged.
}

\begin{document}
\thispagestyle{empty}
{
\title{\thetitle}
\theauthors

\maketitle
\begin{abstract}
{\footnotesize \theabstract}
\end{abstract}

\ifAMS
    {\par\noindent\emph{AMS Subject Classification:} Primary \kwdp. Secondary \kwds}
    {\par\noindent\emph{JEL codes}: \kwdp}

} 

\tableofcontents

\Section{Introduction}
\label{Sdevexpintr}
Let \( \xiv \) be a zero mean Gaussian vector in \( \R^{\dimp} \) for \( \dimp \) large and let
\( \QP \colon \R^{\dimp} \to \R^{\dimq} \) be a linear mapping.
Then \( \QP \xiv \) is also zero mean Gaussian with the variance
\( \BBH = \Var(\QP \xiv) = \QP \Var(\xiv) \QP^{\T} \). 
For the squared norm \( \| \QP \xiv \|^{2} \), it holds \( \E \| \QP \xiv \|^{2} = \tr (\BBH) \), 
\( \Var(\| \QP \xiv \|^{2}) = \tr(\BBH^{2}) \), 
and this random variable concentrates around its expectation \( \tr (\BBH) \) in the sense that for any \( \xx > 0 \)
\begin{EQ}[lcl]
	\P\Bigl( \| \QP \xiv \|^{2} - \tr (\BBH) > 2 \sqrt{\xx \tr(\BBH^{2})} + 2 \| \BBH \| \xx \Bigr)
	& \leq &
	\ex^{-\xx} ,
	\\
	\P\Bigl( \| \QP \xiv \|^{2} - \tr (\BBH) < - 2 \sqrt{\xx \tr(\BBH^{2})} \Bigr)
	& \leq &
	\ex^{-\xx} ; 
\label{v7hfte535ghewjfyw}
\end{EQ}
see e.g. \cite{laurentmassart2000}.
An extension of these bounds to a sub-gaussian case is discussed in \cite{RuVe2013}.
In the recent years, a number of new results were obtained in this direction.
We refer to \cite{KlZi2019} for an extensive overview and advanced results on Hanson-Wright type concentration inequalities.
However, sub-gaussian behavior of the vector \( \xiv \) can be very restrictive in many applications.
Typical examples are given by weighted sums of Bernoulli, Poisson, exponential random variables or by empirical covariance.
For all such examples, the tails of \( \xiv \) are sub-exponential, but su
This note presents some deviation bounds on the norm \( \| \xiv \| \) for the case when the moment generating function
\( \E \exp \langle \xiv,\uv \rangle \) is well defined on a sufficiently large but bounded set of vectors \( \uv \in \R^{\dimp} \).
The main challenge of the study is that the standard technique based on Markov inequality for \( \| \QP \xiv \|^{2} \)
does not apply because the exponential  moments of \( \| \QP \xiv \|^{2} \) diverge.
We apply a kind of trimming technique combined with pilling device to obtain 
nearly sharp bounds in the form
\begin{EQA}
	\P\bigl( \| \QP \xiv \| > \zqc(\BBH,\xx) \bigr)
	& \leq &
	3 \ex^{-\xx}
\label{idcudc8edkjfvy23wkv8ude}
\end{EQA}
with a phase transition effect: 
for \( \xx \leq \xxc \approx \gmb^{2}/4 \), the quantile function \( \zqc \) is exactly as in the Gaussian case:
\( \zqc^{2}(\BBH,\xx) = \tr (\BBH) + 2 \sqrt{\xx \tr(\BBH^{2})} + 2 \| \BBH \| \xx \).
For \( \xx > \xxc \), the function \( \zqc(\BBH,\xx) \) grows linearly as 
\( \CONST \xx/\gmb \) for an absolute constant \( \CONST \).

The paper is organized as follows.
The main deviation bounds on \( \| \QP \xiv \| \) are collected in Section~\ref{Sdevboundexp}.
Applications to weighted Bernoulli vector sums are given in Section~\ref{SdevBern}.
Sharp deviation bounds for empirical covariance matrix are given in Section~\ref{SFrobnorm}.
Some useful technical facts about Gaussian quadratic forms are collected in the Appendix~\ref{SmomentqfG} and
Appendix~\ref{SdevboundGauss}.


\def\zqx{\zq_{\circ}}

\Section{Deviation bounds under light exponential tails}
\label{Sdevboundexp}

Let \( \xiv \) be a zero mean random vector in \( \R^{\dimp} \) with covariance \( \Var(\xiv) \) and 
let \( \QP \colon \R^{\dimp} \to \R^{\dimq} \) be a linear mapping.
This section presents some deviation bounds on the norm \( \| \QP \xiv \| \) for the case of light exponential tails of \( \xiv \).
Namely, 

\begin{description}
\item[\label{gmref} \( \bb{(\gmb)} \)]
	\textit{for some fixed \( \gmb > 0 \) and some self-adjoint operator \( \HVB^{2} \) in \( \R^{\dimp} \) with} 
	\( \HVB^{2} \geq \Var(\xiv) \), 
\begin{EQA}[c]
    \cdens(\uv)
    \eqdef
    \log \E \exp\bigl( \langle \uv, \HVB^{-1} \xiv \rangle \bigr)
    \le
    \frac{\| \uv \|^{2}}{2} \, ,
    \qquad
    \uv \in \R^{\dimp}, \, \| \uv \| \le \gmb ,
\label{expgamgm}
\end{EQA}
\end{description}

In fact, it is sufficient to assume that 
\begin{EQA}
	\sup_{\| \uv \| \leq \gmb} \E \exp\bigl( \langle \uv, \HVB^{-1} \xiv \rangle \bigr)
	& \leq &
	\CONST .
\label{sgagEexlgHx}
\end{EQA}
The quantity \( \CONST \) can be very large but it is not important.
Indeed, the function \( \cdens(\uv) \) is analytic on the disk \( \| \uv \| \leq \gmb \), 
and condition \eqref{sgagEexlgHx} implies an analog of \eqref{expgamgm}: 
\begin{EQA}
	\cdens(\uv) 
	& \leq &
	\frac{\| \uv \|^{2}}{2} +\frac{\dltwu_{3} \| \uv \|^{3}}{6} 
	\leq 
	\frac{\| \uv \|^{2}}{2} \Bigl( 1 + \frac{\dltwu_{3} \gmn}{3} \Bigr)
	\, ,
	\qquad
	\| \uv \| \leq \gmn \, ,
\label{0hy544rtgte7hjuwghbj}
\end{EQA}
for a fixed value \( \dltwu_{3} \).
Moreover, reducing \( \gmb \) allows to take \( \HVB^{2} \) equal or close to \( \Var(\xiv) \) and \( \dltwu_{3} \)
close to zero.
The next section presents our main results under \nameref{gmref}.
The proofs are postponed until the end of the section.

\Subsection{Main results}
Let a random vector \( \xiv \) satisfy \( \E \xiv = 0 \) and \nameref{gmref}. 
The goal is to establish possibly sharp deviation bounds on \( \| \QP \xiv \|^{2} \)
for a given linear mapping \( \QP \colon \R^{\dimp} \to \R^{\dimq} \).
Define 
\begin{EQ}[rcl]
	\BBH 
	& \eqdef &
	\QP \HVB^{2} \QP^{\T}, \quad
	\dimH 
	\eqdef  
	\tr(\BBH) , 
	\quad
	\vH^{2} 
	\eqdef
	\tr(\BBH^{2}) , 
	\quad 
	\supH 
	\eqdef 
	\| \BBH \| ,
	\\
	\zq^{2}(\BBH,\xx) 
	& \eqdef & 
	\tr \BBH + 2 \sqrt{\xx \tr(\BBH^{2})} + 2 \xx \| \BBH \|
	=
	\dimH + 2 \vH \sqrt{\xx} + 2 \xx \supH .
\label{h5ete5f5tegvhvcdvhyendf}
\end{EQ} 
Also fix some \( \rhoH < 1 \), a standard choice is \( \rhoH = 1/2 \).
Our main result applies for all \( \xx \) satisfying the condition
\begin{EQA}
	\zq^{2}(\BBH,\xx)
	& \leq &
	\rhoH \, \biggl(\frac{\gmb \sqrt{\supH}}{\muH(\xx)} - \sqrt{\frac{\dimH}{\muH(\xx)}} \biggr)^{2}
\label{kv7367ehjgruwwcewyd}
\end{EQA}
with \( \zq(\BBH,\xx) \) from \eqref{h5ete5f5tegvhvcdvhyendf} and
\( \muH(\xx) \) defined by \( \muH^{-1}(\xx) = 1 + \frac{\vH}{2 \supH \sqrt{\xx}} \);
see \eqref{1v2sxm12m1m}.
One can see that the left hand-side of \eqref{kv7367ehjgruwwcewyd} increases with \( \xx \) while 
the right hand-side decreases.
Therefore, there exists a unique root \( \xxc \) such that with \( \muHc = \muH(\xxc) \)
\begin{EQA}
	\zq^{2}(\BBH,\xxc)
	& = &
	\rhoH \, \biggl(\frac{\gmb \sqrt{\supH}}{\muHc} - \sqrt{\frac{\dimH}{\muHc}} \biggr)^{2} .
\label{kv7367ehjgruwwcewyde}
\end{EQA}
The value \( \xxc \) is important, it describes the \emph{phase transition} effect:
the upper quantile function of \( \| \QP \xiv \| \) exhibits the Gaussian-like behavior for \( \xx \leq  \xxc \), 
while it grows linearly with \( \xx/\gmb \) for \( \xx > \xxc \) as in a sub-exponential case.


\begin{theorem}
\label{Tdevboundgm}
Assume \nameref{gmref}.
Fix \( \xxc \) by \eqref{kv7367ehjgruwwcewyde} for some 
\( \rhoH \leq 1/2 \).
It holds 
\begin{EQA}
    \P\bigl( \| \QP \xiv \| \ge \zq(\BBH,\xx) \bigr)
    & \le &
    3 \ex^{-\xx} ,
    \qquad
    \xx \leq \xxc \, .
\label{PxivbzzBBroB}
\end{EQA}    
For \( \rhoH = 1/2 \), the value \( \xxc \) from \eqref{kv7367ehjgruwwcewyde} and \eqref{PxivbzzBBroB} fulfills 
\begin{EQA}
	\frac{1}{4} \biggl( \gmb - \sqrt{\frac{2\dimH}{\supH}} \biggr)_{+}^{2}
	\leq 
	\xxc 
	& \leq &
	\frac{\gmb^{2}}{4} \, .
\label{if7h3rhgy4676rfhdsjw}
\end{EQA}
If \( \gmb > \sqrt{2 \dimH / \supH} \) then \( \zqc = \zq(\BBH,\xxc) \) follows
\begin{EQA}
	\gmb \sqrt{\supH /2} - (1 - 2^{-1/2}) \sqrt{\dimH} 
	\leq 
	\zqc
	& \leq &
	\gmb \sqrt{\supH /2} + \sqrt{\dimH} \, .
\label{f9oi4eogjdgtvgyuj4ek}
\end{EQA}
\end{theorem}


The results of Theorem~\ref{Tdevboundgm} state nearly Gaussian deviation bounds for 
the norm of the vector \( \QP \xiv \) satisfying \nameref{gmref}. 
Namely, the Gaussian deviation bound \( \P\bigl( \| \QP \xiv \| \ge \zq(\BBH,\xx) \bigr) \le \ex^{-\xx} \)
from Theorem~\ref{TexpbLGA} applies with the additional factor 3 for all \( \xx \leq \xxc \).
Condition \( \gmb \gg \sqrt{\dimH / \supH} \) is important.
Otherwise, the value \( \xxc \) is not significantly large and the zone \( \xx \leq \xxc \) 
with Gaussian-like quantiles is too narrow.
It turns out that out of this range, the norm \( \| \QP \xiv \| \) exhibits a sub-exponential behavior.

\begin{theorem}
\label{TQPxivlarge}
Assume \nameref{gmref}.
With \( \xxc \) from \eqref{kv7367ehjgruwwcewyde} and \( \zqc = \zq(\BBH,\xxc) \),
set \( \constg = \frac{\sqrt{\rhoH} \, \gmb}{(2 + \sqrt{\rhoH}) \sqrt{\supH}} \).
It holds 
\begin{EQ}[lll]
	&
	\P\bigl( \| \QP \xiv \| > \zqc + \constg^{-1} (\xx - \xxc) \bigr)
	\leq 
	3 \ex^{-\xx} ,
	&
	\quad
	\xx \geq \xxc \, ,
	\\
	&
	\P\bigl( \| \QP \xiv \| > \zq \bigr)
	\leq
	3 \exp\{ -\xxc - \constg (\zq - \zqc) \} ,
	&
	\quad
	\zq \geq \zqc \, .
\label{9vkjv6njbih9t69t3wfg}
\end{EQ}
\end{theorem}

The obtained deviation bounds of Theorem~\ref{Tdevboundgm} and Theorem~\ref{TQPxivlarge}
can be fused into one.
To be more specific, we fix \( \rhoH = 1/2 \).

\begin{corollary}
\label{CTQPxivlarge}
Assume \nameref{gmref}. 
Let \( \xxc \) be defined by \eqref{kv7367ehjgruwwcewyde} with \( \rhoH = 1/2 \).
For all \( \xx > 0 \)
\begin{EQA}
	\P\bigl( \| \QP \xiv \| > \zqc(\BBH,\xx) \bigr)
	& \leq &
	3 \ex^{-\xx} ,
\label{uyfuyerd7e7uhh8yy689t}
\end{EQA}
where with \( \constg \eqdef \frac{\gmb}{(\sqrt{8} + 1) \sqrt{\supH}} \)
and \( \xx \wedge \xxc \eqdef \min\{ \xx , \xxc \} \)
\begin{EQA}
	\zqc(\BBH,\xx)
	& \eqdef &
	\zq(\BBH,\xx \wedge \xxc)
	+ \constg^{-1} (\xx - \xxc)_{+}
	=
	\begin{cases}
		\zq(\BBH,\xx), & \xx \leq \xxc \, ,
		\\
		\zq(\BBH,\xxc) + \dfrac{\xx - \xxc}{\constg} , & \xx > \xxc \, .
	\end{cases}
	\qquad
	\qquad
\label{dyv6ejf8gjwkerih83}
\end{EQA}
Moreover, \( \xxc \) follows \eqref{if7h3rhgy4676rfhdsjw} and \( \zqc = \zq(\BBH,\xxc) \) 
satisfies \eqref{f9oi4eogjdgtvgyuj4ek} provided \( \gmb \geq \sqrt{2\dimH / \supH} \).
\end{corollary}

If \( \gmb \gg \sqrt{\dimH / \supH} \) then \( \xxc \) is large and 
\( \zqc(\BBH,\xx) = \zq(\BBH,\xx) \leq \sqrt{\dimH} + \sqrt{2 \xx \supH} \) for all reasonable \( \xx \).
For \( \gmb < \sqrt{2 \dimH / \supH} \), the accurate bound \eqref{dyv6ejf8gjwkerih83} can be simplified by a linear majorant
which does not involve \( \xxc \).

\begin{theorem}
\label{Llin3maj}
Assume \nameref{gmref}. 
Fix \( \constg = \frac{\gmb }{(\sqrt{8} + 1) \sqrt{\supH}} \).
Then \eqref{uyfuyerd7e7uhh8yy689t} applies with
\begin{EQA}
	\zqc(\BBH,\xx)
	& \leq &
	\sqrt{\dimH} + \frac{\constg}{\sqrt{2}} + \constg^{-1} \xx \, .
\label{uv7ey4hb743h3hiti4ffd}
\end{EQA}
\end{theorem}

The next result provides some upper bounds on the exponential moments of \( \| \QP \xiv \| \).
We distinguish between zones \( \zq \leq \zqc \) and \( \zq > \zqc \) with \( \zqc = \zq(\BBH,\xxc) \);
see \eqref{kv7367ehjgruwwcewyde}.

\begin{theorem}
\label{TQPlargex}
Assume \nameref{gmref}.
Let \( \xxc \) fulfill \eqref{kv7367ehjgruwwcewyde} and \( \zqc = \zq(\BBH,\xxc) \).
For any \( \zq \in [\sqrt{\dimH}, \zqc] \) and any \( \nuH \leq \frac{\zq - \sqrt{\dimH}}{2\sqrt{\supH}} \), it holds
\begin{EQA}
	\E \ex^{\nuH \| \QP \xiv \|} \Ind(\| \QP \xiv \| \geq \zq)
	& \leq &
	6 \exp\Bigl\{  \nuH \zq - \frac{(\zq - \sqrt{\dimH})^{2}}{2 \supH} \Bigr\} .
\label{ufikwk3ei9vgkj4k4giw3wl}
\end{EQA}
Further, for any \( \nuH < \constg \eqdef \frac{\gmb \sqrt{\rhoH}}{\sqrt{\supH} \, (2 + \sqrt{\rhoH})} \)
\begin{EQA}
	\E \ex^{\nuH \| \QP \xiv \| } \Ind(\| \QP \xiv \| > \zqc)
	& \leq &
	\frac{3 \constg}{\constg - \nuH} \, 
	\exp\biggl\{  \nuH \zqc - \frac{(\zqc - \sqrt{\dimH})^{2}}{2 \supH} \biggr\} \, .
\label{jhf7yehruybyrhe3wevire}
\end{EQA}
Moreover, for \( \zq \geq \zqc \)
\begin{EQA}
	\E \ex^{\nuH \| \QP \xiv \| } \Ind(\| \QP \xiv \| > \zq)
	& \leq &
	\frac{3 \constg}{\constg - \nuH} \, 
	\exp \Bigl\{ \nuH \zqc - \frac{(\zqc - \sqrt{\dimH})^{2}}{2 \supH} - (\constg - \nuH) (\zq - \zqc) \Bigr\} \, .
	\qquad
\label{jhf7yehruybyrhe3wevire2}
\end{EQA}
\end{theorem}

\Subsection{Proof of Theorem~\ref{Tdevboundgm}}
By normalization, one can easily reduce the study to the case \( \| \BBH \| = 1 \).
Moreover, replacing \( \xiv \) with \( \HVB^{-1} \xiv \) and \( \QP \) with \( \QP \HVB \) reduces the proof to the situation 
with \( \HVB = \Id_{\dimp} \).
This will be assumed later on.
%
For \( \muH \in (0,1) \) and \( \zzH(\muH) = \gmb / \muH - \sqrt{\dimH/\muH} > 0 \),
define trimming \( t_{\muH}(\uv) \) of \( \uv \in \R^{\dimp} \) as
\begin{EQA}
	t_{\muH}(\uv)
	& \eqdef &
	\begin{cases}
		\uv, & \text{ if } \| \uv \| \leq \zzH(\muH), 
		\\
		\frac{\zzH(\muH)}{\| \uv \|} \, \uv, & \text{ otherwise}. 
	\end{cases}
\label{ojchcy63eyfvey3gcbfkg}
\end{EQA}
By construction \( \| t_{\muH}(\uv) \| \leq \zzH(\muH) \) for all \( \uv \in \R^{\dimp} \). 

\begin{lemma}
\label{LGDBqfexpB}
Assume \nameref{gmref} and let \( \| \BBH \| = 1 \).
Fix \( \muH \in (0,1) \) s.t. \( \zzH(\muH) = \gmb / \muH - \sqrt{\dimH/\muH} > 0 \).
Then with \( t_{\muH}(\cdot) \) from \eqref{ojchcy63eyfvey3gcbfkg} 
\begin{EQA}
	\E \exp\Bigl\{ \frac{\muH}{2} \, t_{\muH}^{2}(\QP \xiv) \Bigr\}
	& \leq &
	2 \exp\{ \Pmuvp(\muH) \} ,
	\qquad
\label{wBmu2vA241mmutrB}
\end{EQA}
where
\begin{EQA}
	\Pmuvp(\muH)
	& \eqdef &
	\frac{\muH^{2} \vH^{2}}{4 (1 - \muH)} + \frac{\muH \, \dimH}{2} \, .
\label{sdfw6rewqrwdwtffwet66}
\end{EQA}
Furthermore, for any \( \zzH < \zzH(\muH) \)
\begin{EQA}
	\P\bigl( \| \QP \xiv \| > \zzH, \| \QP \xiv \| \leq \zzH(\muH) \bigr)
	& \leq &
	2 \exp\Bigl\{ - \frac{\muH \, \zzH^{2}}{2} + \Pmuvp(\muH) \Bigr\} .
\label{9cdkcf736eryghj7y34wwde}
\end{EQA}
\end{lemma}

\begin{proof}
Let us fix any value of \( \xiv \).
We intend to show that
\begin{EQA}
	\exp\Bigl\{ \frac{\muH}{2} \, \| t_{\muH}(\QP \xiv) \|^{2} \Bigr\}
	& \leq &
	2 \Egs \, \exp\{ \muH^{1/2} \gaussv^{\T} t_{\muH}(\QP \xiv) \} .
\label{0xcujfc5nfiuf76w3hfc}
\end{EQA}
Here \( \Egs \) means conditional expectation w.r.t. \( \gaussv \sim \ND(0,\Id_{\dimp}) \) given \( \xiv \).
Obviously, with \( A = \{ \uv \colon \muH^{1/2} \| \QP^{\T} \uv \| \leq \gmb \} \), 
it suffices to check that
\begin{EQA}
	\II_{\muH}(\xiv)
	& \eqdef &
	\Egs \exp \Bigl\{ \muH^{1/2} \gaussv^{\T} t_{\muH}(\QP \xiv) - \frac{\muH}{2} \| t_{\muH}(\QP \xiv) \|^{2} \Bigr\} 
	\Ind(\gaussv \in A)
	\geq 
	1/2 .
\label{ggfiokhyhiu8y6eeeewe}
\end{EQA}
With \( \CONSTi_{\dimp} = (2\pi)^{-\dimp/2} \), it holds
\begin{EQA}
	\II_{\muH}(\xiv)
	&=&
	\CONSTi_{\dimp} \int_{A} \exp \Bigl( 
		\muH^{1/2} \uv^{\T} t_{\muH}(\QP \xiv) 
		- \frac{\muH}{2} \| t_{\muH}(\QP \xiv) \|^{2} - \frac{1}{2} \| \uv \|^{2} 
	\Bigr) d\uv
	\\
	&=&
	\CONSTi_{\dimp} \int_{A} \exp \Bigl( - \frac{1}{2} \| \uv - \muH^{1/2} t_{\muH}(\QP \xiv) \|^{2} \Bigr) d\uv
	=
	\Pgs(\gaussv - \muH^{1/2} t_{\muH}(\QP \xiv) \in A ) .
\label{kj8vfuff7f7f7ywwew}
\end{EQA}
The definition of \( A \) and the condition \( \| t_{\muH}(\QP \xiv) \| \leq \zzH(\muH) \) imply
in view of \( \| \QP \| \leq 1 \)
\begin{EQA}
	&& \nquad
	\Pgs(\gaussv - \muH^{1/2} t_{\muH}(\QP \xiv) \in A )
	=
	\Pgs\bigl( \| \QP^{\T} (\gaussv - \muH^{1/2} t_{\muH}(\QP \xiv)) \| \leq \gmb/\muH^{1/2} \bigr)
	\\
	& \geq &
	\Pgs\bigl( \| \QP^{\T} \gaussv \| \leq \gmb/\muH^{1/2} - \muH^{1/2} \zzH(\muH) \bigr)
	\geq 
	\Pgs\bigl( \| \QP^{\T} \gaussv \| \leq \sqrt{\dimH} \bigr)
	\geq 
	1/2
\label{mnriuvu8frt7r77r7urfu}
\end{EQA}
and \eqref{ggfiokhyhiu8y6eeeewe} follows.
Taking expectation for both sides of \eqref{0xcujfc5nfiuf76w3hfc} and the use of Fubini's theorem yield
\begin{EQA}
	&& \nquad
	\E \exp\Bigl\{ \frac{\muH}{2} \, \| t_{\muH}(\QP \xiv) \|^{2} \Bigr\} 
	\leq 
	2 \Egs \bigl\{ \E \exp\{ \muH^{1/2} \gaussv^{\T} t_{\muH}(\QP \xiv) \} 
	\, \Ind(\muH^{1/2} \| \QP^{\T} \gaussv \| \leq \gmb) \bigr\} .
\label{icjfcufry7fvc6e32w33wr}
\end{EQA}
Obviously, for any \( \uv \in \R^{\dimp} \)
\begin{EQA}
	\exp\{ \uv^{\T} t_{\muH}(\QP \xiv) \} + \exp\{ - \uv^{\T} t_{\muH}(\QP \xiv) \} 
	& \leq &
	\exp\{ \uv^{\T} \QP \xiv \} + \exp\{ - \uv^{\T} \QP \xiv \}
\label{odfu87e76376rhfyt6wh}
\end{EQA} 
for any \( \uv \in \R^{\dimp} \) with \( \| \QP^{\T} \uv \| \leq \gmb \)
and by \eqref{expgamgm} 
\begin{EQA}
	\E \exp\Bigl\{ \frac{\muH}{2} \, \| t_{\muH}(\QP \xiv) \|^{2} \Bigr\}
	& \leq &
	2 \Egs \Bigl\{ \exp\Bigl( \frac{1}{2} \, \| \muH^{1/2} \gaussv^{\T} \QP \|^{2} \Bigr) \, 
		\Ind(\muH^{1/2} \| \QP^{\T} \gaussv \| \leq \gmb) \Bigr\}
	\\
	& \leq &
	2 \Egs \, \exp\Bigl( \frac{1}{2} \, \| \muH^{1/2} \gaussv^{\T} \QP \|^{2} \Bigr)
	=
	2 \det(\Id_{\dimp} - \muH \QP^{\T} \QP)^{-1/2} .
\label{8ji89656hk0wchdne83}
\end{EQA}
We also use that for any \( \muH > 0 \) by \eqref{m2v241m4mj1p},
\begin{EQA}
	\log \det\bigl( \Id - \muH \BBH \bigr)^{-1/2} 
	& \leq &
	\frac{\muH \tr (\BBH)}{2} + \frac{\muH^{2} \tr (\BBH^{2})}{4 (1 - \muH)} 
	=
	\Pmuvp(\muH) \, ,
\label{mu2v241mmiulogIm12}
\end{EQA}
and the first statement follows.
Moreover, by Markov's inequality 
\begin{EQA}
	\P\bigl( \| \QP \xiv \| > \zzH, \| \QP \xiv \| \leq \zzH(\muH) \bigr)
	& \leq &
	\ex^{- \muH \, \zzH^{2}/2} \E \exp\Bigl\{ \frac{\muH}{2} \, \| t_{\muH}(\QP \xiv) \|^{2} \Bigr\} 
	\leq 
	2 \exp\Bigl\{ - \frac{\muH \, \zzH^{2}}{2} + \Pmuvp(\muH) \Bigr\} ,
\label{9cdkcf736eryfvkler4895ty34wwde}
\end{EQA}
and \eqref{9cdkcf736eryghj7y34wwde} follows as well.
\end{proof}

The use of \( \muH = \muH(\xx) \) from \eqref{1v2sxm12m1m} in \eqref{wBmu2vA241mmutrB} yields 
\begin{EQA}
	- \frac{\muH \zq^{2}(\BBH,\xx)}{2} + \Pmuvp(\muH)
	&=&
	- \xx \, ,
\label{uv73hf8er74e7eereew}
\end{EQA}
and similarly to the proof of Theorem~\ref{TexpbLGA}
\begin{EQA}
	\P\Bigl( \| \QP \xiv \|^{2} > \zq^{2}(\BBH,\xx), \,
		\| \QP \xiv \| \leq \zzH(\muH)
	 \Bigr)
	 & \leq &
	 2 \ex^{-\xx} .
\label{2emxPblrHm1B}
\end{EQA}
It remains to consider the probability of large deviation
\( \P\bigl( \| \QP \xiv \| > \zzH(\muH) \bigr) \).

\begin{lemma}
\label{Ldvbetagmb}
Assume \( \| \BBH \| = 1 \).
Given \( \xx > 0 \), fix \( \muH = \muH(\xx) \) and \( \zzH(\muH) = \gmb/\muH - \sqrt{\dimH/\muH} \).
Assume \eqref{kv7367ehjgruwwcewyd} for some \( \rhoH \leq 1/2 \).
Then 
\begin{EQA}
	\P\bigl( \| \QP \xiv \| > \zzH(\muH) \bigr)
	& \leq &
	\ex^{-\xx} .	
\label{yvy3n3eubvuj2tcy3he}
\end{EQA}
\end{lemma}

\begin{proof}
Denote \( \eta = \| \QP \xiv \| \).
By \eqref{2emxPblrHm1B}
\begin{EQA}
	\P\Bigl( \eta > \zq(\BBH,\xx), \,
		\eta \leq \zzH(\muH)
	 \Bigr)
	 & \leq &
	 2 \ex^{-\xx} ,
\label{me22Iezmemcz2}
\end{EQA}
For \( \muH = \muH(\xx) \), it holds \eqref{uv73hf8er74e7eereew} with \( \Pmuvp(\muH) \) 
given by \eqref{sdfw6rewqrwdwtffwet66}.
Bounding the tails of \( \eta \) in the region \( \eta > \zzH(\muH) \) requires another choice of \( \muH \).
Namely, we apply \eqref{9cdkcf736eryghj7y34wwde} with \( \rhoH \muH \) instead of \( \muH \) yielding
\begin{EQA}
	\P\bigl( \eta > \zzH(\muH), \eta \leq \zz(\rhoH \muH) \bigr)
	& \leq &
	2 \exp\Bigl\{ - \frac{\rhoH \muH \, \zz^{2}(\muH)}{2} + \Pmuvp(\rhoH \muH) \Bigr\} .
\label{uv7e3j2kv88eu3e3536}
\end{EQA}
In a similar way, applying \eqref{me22Iezmemcz2} with \( \rhoH^{2} \muH \) in place of \( \muH \)
and using that 
\begin{EQA}
	\rhoH \, \zz(\rhoH \muH)
	&=&
	\gmb/\muH - \sqrt{\rhoH \, \dimH/\muH}
	\leq 
	\zzH(\muH) 
\label{f9i2mg9gj3g2553}
\end{EQA}
yields
\begin{EQA}
	&& \nquad
	\P\bigl( \eta > \zz(\rhoH \muH), \eta \leq \zz(\rhoH^{2} \muH) \bigr)
	\leq 
	2 \exp\Bigl\{ - \frac{\rhoH^{2} \muH \, \zz^{2}(\rhoH \muH)}{2} + \Pmuvp(\rhoH^{2} \muH) \Bigr\} 
	\\
	& \leq &
	2 \exp\Bigl\{ - \frac{\muH \, \zz^{2}(\muH)}{2} + \Pmuvp(\rhoH^{2} \muH) \Bigr\} .
\label{t6dxhywywybcxcttw2q}
\end{EQA}
This trick can be applied again and again yielding in view of \eqref{f9i2mg9gj3g2553}
\begin{EQA}
	\P\bigl( \eta > \zzH(\muH) \bigr)
	& \leq &
	\sum_{k=0}^{\infty} \P\bigl( \eta > \zz(\rhoH^{k} \muH), \eta \leq \zz(\rhoH^{k+1} \muH) \bigr)
	\\
	& \leq &
	\sum_{k=0}^{\infty} 2 \exp\bigl\{ - \rhoH^{k+1} \muH \, \zz^{2}(\rhoH^{k} \muH)/2 + \Pmuvp(\rhoH^{k+1} \muH) \bigr\} 
	\\
	& \leq &
	\sum_{k=0}^{\infty}  2 \exp\bigl\{ - \rhoH^{-k+1} \muH \, \zz^{2}(\muH)/2 + \Pmuvp(\rhoH^{k+1} \muH) \bigr\} .
\label{uenruvueju33j2223figgi}
\end{EQA}
Condition \( \rhoH \, \zz^{2}(\muH) \geq \zq^{2}(\BBH,\muH)/2 \) and \eqref{uv73hf8er74e7eereew} ensure
for \( \rhoH \leq 1/2 \)
\begin{EQA}
	\P\bigl( \eta > \zzH(\muH) \bigr)
	& \leq &
	\sum_{k=0}^{\infty} 2 \exp\bigl\{ - \rhoH^{-k} \muH \, \zq^{2}(\BBH,\muH)/2 + \Pmuvp(\rhoH^{k+1} \muH) \bigr\}
	\\
	& \leq &	
	2 \sum_{k=0}^{\infty}  \exp\bigl\{ \Pmuvp(\rhoH^{k+1} \muH) - \rhoH^{-k} \Pmuvp(\muH) - \rhoH^{-k} \xx \bigr\} 
	\leq 
	\ex^{-\xx} .	
\label{yvy3n3eubvuj2tyche}
\end{EQA}
This yields \eqref{yvy3n3eubvuj2tcy3he}.
\end{proof}
Putting together \eqref{2emxPblrHm1B} and \eqref{yvy3n3eubvuj2tcy3he} yields \eqref{PxivbzzBBroB}.

Now we check \eqref{if7h3rhgy4676rfhdsjw}.
Normalization by \( \supH \) reduces the proof to the case \( \| \BBH \| = \| \QP \HVB^{2} \QP^{\T} \| = 1 \).
We use the simplified bounds 
\( \zq(\BBH,\xx) \leq \sqrt{\dimH} + \sqrt{2\xx} \) and \( \muH^{-1} = 1 + \sqrt{\dimH/(4\xx)} \). 
Now \eqref{kv7367ehjgruwwcewyd} with \( \rhoH = 1/2 \) can be rewritten as
\begin{EQA}
	\gmb
	& \geq &
	\sqrt{\muH \, \dimH} + \muH \sqrt{2} \bigl( \sqrt{\dimH} + \sqrt{2\xx} \bigr) .
\label{hdf6eh3ye4545dwwe}
\end{EQA}
The use of \( \muH = \sqrt{4 \xx} / (\sqrt{4\xx} + \sqrt{\dimH}) \) yields
\begin{EQA}
	\muH \sqrt{2} \bigl( \sqrt{\dimH} + \sqrt{2\xx} \bigr)
	&=&
	\sqrt{8\xx} \,\, \frac{\sqrt{\dimH} + \sqrt{2\xx}}{\sqrt{\dimH} + \sqrt{4\xx}}
	\geq \sqrt{4\xx} \, ,
\label{uviuedu8e347y4rt87hiuy}
\end{EQA}
and \eqref{hdf6eh3ye4545dwwe} is not possible for \( \xx > \gmb^{2}/4 \).
Further, with \( \yy = \sqrt{4\xx}/\gmb \) and \( \alp = \sqrt{\dimH} / \gmb \)
\begin{EQA}
	\frac{\sqrt{\muH \, \dimH} + \muH \sqrt{2} \bigl( \sqrt{\dimH} + \sqrt{2\xx} \bigr)}{\gmb}
	&=&
	\sqrt{\frac{\yy \alp^{2}}{\alp+\yy}} + \frac{\yy (\sqrt{2} \alp + \yy)}{\alp+\yy}
	\leq 
	\alp + \yy + \frac{\yy(\sqrt{2} - 1) \alp}{\alp+\yy} 
	\leq 
	\yy + \sqrt{2} \alp .
\label{jhfy7637ur8thj986eewe}
\end{EQA}
Together with \eqref{hdf6eh3ye4545dwwe}, this yields \( \yy \geq 1 - \sqrt{2} \alp \) and \eqref{if7h3rhgy4676rfhdsjw} follows.
For \eqref{f9oi4eogjdgtvgyuj4ek} we use \( \zqc \leq \sqrt{\dimH} + \sqrt{2 \supH \xxc} \)
and \( \zqc \geq \sqrt{\dimH/2} + \sqrt{2 \supH \xxc} \).

\Subsection{Proof of Theorem~\ref{TQPxivlarge}}
Assume w.l.o.g. \( \supH = 1 \).
First we present an accurate deviation bound, which, however, does not provide a closed form quantile function 
for \( \| \QP \xiv \| \).
Then we show how it implies a rough linear upper bound on this quantile function.
For \( \xxc \) from \eqref{kv7367ehjgruwwcewyde} and \( \xx > \xxc \), fix \( \muH \) by the relation   
\begin{EQA}
	\frac{\rhoH \, \muH \, \zzH^{2}(\muH)}{2} 
	&=& 
	\xx + \Pmuvp(\muH) 
	=
	\xx + \frac{\muH \, \dimH}{2} + \frac{\muH^{2} \vH^{2}}{4 (1 - \muH)} \, ,
\label{jejfye4ye3hfnw3jfeu}
\end{EQA}
where \( \zzH(\muH) = \gmb/\muH - \sqrt{\dimH/\muH} \)%
; cf. \eqref{uv73hf8er74e7eereew}.
It is easy to see that the solution \( \muH \) exists and unique.
Moreover, if \( \xx = \xxc \) then \( \muH = \muHc \) and \( \zzH^{2}(\muHc) = \zq^{2}(\BBH,\xxc) \);
see \eqref{kv7367ehjgruwwcewyde}.
If \( \xx > \xxc \), then \( \muH < \muHc \) and \( \zzH^{2}(\muH) > \zq^{2}(\BBH,\xx) \).

\begin{lemma}
\label{Ldevbosub}
For \( \xx > \xxc \), define \( \muH \) by \eqref{jejfye4ye3hfnw3jfeu}.
Then with \( \zzH(\muH) = \gmb/ \muH - \sqrt{\dimH/\muH} \)
\begin{EQA}
	\P\bigl( \| \QP \xiv \|^{2} > \rhoH \, \zzH^{2}(\muH) \bigr)
	& \leq &
	3 \ex^{-\xx} \, .
\label{uvw7erf38sjcwtet23w3}
\end{EQA}
\end{lemma}
\begin{proof}
We again apply Lemma~\ref{LGDBqfexpB}, however, the choice \( \muH = \muH(\xx) \) from \eqref{1v2sxm12m1m} 
is not possible anymore in view of \( \zq(\BBH,\xx) > \zzH(\muH) \). 
More precisely, for \( \xx \) large, the value \( \muH(\xx) \) approaches one and this choice of \( \muH \) 
yields the value \( \zzH(\muH) \) smaller than we need.
To cope with this problem, we apply \eqref{9cdkcf736eryghj7y34wwde} of Lemma~\ref{LGDBqfexpB} with a sub-optimal \( \muH \) 
from \eqref{jejfye4ye3hfnw3jfeu} ensuring 
\( \rhoH \muH \, \zzH^{2}(\mu) - \Pmuvp(\muH) = \xx \).
By \eqref{9cdkcf736eryghj7y34wwde} of Lemma~\ref{LGDBqfexpB} 
\begin{EQA}
	\P\bigl( \| \QP \xiv \| > \sqrt{\rhoH} \, \zzH(\muH), \| \QP \xiv \| \leq \zzH(\muH) \bigr)
	& \leq &
	2 \exp\Bigl\{ - \frac{\rhoH \muH \, \zzH^{2}(\muH)}{2} + \Pmuvp(\muH) \bigr\}
	=
	2 \ex^{-\xx} .
\label{mcyer6e46rf62wuvyrddde}
\end{EQA}
Repeating the arguments from the proof of Lemma~\ref{Ldvbetagmb} implies
\begin{EQA}
	&& \nquad
	\P\bigl( \| \QP \xiv \|^{2} > \rhoH \, \zzH^{2}(\muH) \bigr)
	\leq 
	\sum_{k=0}^{\infty} 2 \exp\Bigl\{ - \frac{1}{2} \rhoH^{k+1} \muH \, \zz^{2}(\rhoH^{k} \muH) + \Pmuvp(\rhoH^{k} \muH) \Bigr\} 
	\\
	& \leq &
	\sum_{k=0}^{\infty}  2 \exp\Bigl\{ - \frac{1}{2} \rhoH^{-k+1} \muH \, \zz^{2}(\muH) + \Pmuvp(\rhoH^{k} \muH) \Bigr\} 
	\\
	& \leq &
	2 \ex^{-\xx} 
	+ 2 \ex^{-\xx} \sum_{k=1}^{\infty}  
		\exp\Bigl\{ - \frac{1}{2} (\rhoH^{-k} - 1) \rhoH \muH \, \zz^{2}(\muH) + \Pmuvp(\rhoH^{k} \muH) - \Pmuvp(\muH) \Bigr\} 
	\leq 
	3 \ex^{-\xx} .
\label{uenruvueju3kH223figgi}
\end{EQA}
as stated in \eqref{uvw7erf38sjcwtet23w3}.
\end{proof}

\noindent
It remains to evaluate \( \rhoH \, \zzH^{2}(\muH) \) with \( \muH \) from \eqref{jejfye4ye3hfnw3jfeu} and
\( \zzH(\muH) = \gmb/\muH - \sqrt{\dimH/\muH} \).
For \( \muH \leq \muHc \)
\begin{EQA}
	\frac{\rhoH}{2} \biggl( \frac{\gmb}{\sqrt{\muH}} - \sqrt{\dimH} \biggr)^{2}
	&=&
	\xx + \Pmuvp(\muH)
\label{uedu7euyewhfvhfrt4qww}
\end{EQA}
and 
\begin{EQA}
	\frac{\sqrt{\rhoH} \, \gmb}{\sqrt{\muH}}
	& = &
	\sqrt{2\xx + 2\Pmuvp(\muH)} + \sqrt{\rhoH \, \dimH} 
	\, .
\label{iujk37bujreyvjwtcjeww}
\end{EQA}
This results in
\begin{EQA}
	\sqrt{\rhoH} \, \zzH(\muH)
	& = &
	\frac{\sqrt{\rhoH}}{\sqrt{\muH}} \biggl( \frac{\gmb}{\sqrt{\muH}} - \sqrt{\dimH} \biggr)
	\leq 
	\frac{1}{\sqrt{\rhoH} \, \gmb} \bigl( \sqrt{2\xx + 2\Pmuvp(\muH)} + \sqrt{\rhoH \, \dimH} \bigr) \, \sqrt{2\xx + 2\Pmuvp(\muH)} 
	\\
	& \leq &
	\frac{1}{\sqrt{\rhoH} \, \gmb} 
	\bigl( 2\xx + 2\Pmuvp(\muHc) + \sqrt{\rhoH \, \dimH(2\xx + 2\Pmuvp(\muHc))} \bigr) 
	\eqdef
	\zqs(\xx)
	\, .
\label{hv7jwwsuvjqs89c83jejfy}
\end{EQA}
By \eqref{kv7367ehjgruwwcewyde}, this inequality becomes equality for \( \xx = \xxc \) and \( \muH = \muHc \) with 
\( \sqrt{\rhoH} \, \zzH(\muHc) = \zqs(\xxc) = \zq(\BBH,\xxc) \).
Furthermore, the derivative of \( \zqs(\xx) \) w.r.t. \( \xx \) satisfies
\begin{EQA}
	\frac{d}{d\xx} \, \zqs(\xx)
	&=&
	\frac{1}{\sqrt{\rhoH} \, \gmb} \biggl( 2 + \frac{\sqrt{\rhoH \, \dimH}}{\sqrt{2\xx + 2\Pmuvp(\muHc)}} \biggr)
	\leq 
	\frac{1}{\sqrt{\rhoH} \, \gmb} \biggl( 2 + \frac{\sqrt{\rhoH \, \dimH}}{\sqrt{2\xxc + 2\Pmuvp(\muHc)}} \biggr) \, .
\label{igk373jfgho6952wq2q2}
\end{EQA}
Moreover, \( 2\xxc + 2\Pmuvp(\muHc) = \zq^{2}(\BBH,\xxc) \) and
\begin{EQA}
	\frac{d}{d\xx} \, \zqs(\xx)
	& \leq &
	\frac{1}{\sqrt{\rhoH} \, \gmb} \biggl( 2 + \frac{\sqrt{\rhoH \, \dimH}}{\zq(\BBH,\xxc)} \biggr)
	\leq 
	\frac{2 + \sqrt{\rhoH}}{\sqrt{\rhoH} \, \gmb} 
\label{h8jkew3fbo9rfkw3tehw32h}
\end{EQA}
yielding
\begin{EQA}
	\zqs(\xx)
	& \leq &
	\zqs(\xxc) + \frac{2 + \sqrt{\rhoH}}{\sqrt{\rhoH} \, \gmb} (\xx - \xxc)
	=
	\zq(\BBH,\xxc) + \frac{2 + \sqrt{\rhoH}}{\sqrt{\rhoH} \, \gmb} (\xx - \xxc)
\label{ujw289fjrej3gyhr2wmdx}
\end{EQA}
and hence, 
\begin{EQA}
	\sqrt{\rhoH} \, \zzH(\muH)
	& \leq &
	\zq(\BBH,\xxc) + \frac{2 + \sqrt{\rhoH}}{\sqrt{\rhoH} \, \gmb} (\xx - \xxc) 
	=
	\zqc + \frac{\xx - \xxc}{\constg} \, .
\label{9fj36vhje37gjgr33erfg}
\end{EQA}
This implies \eqref{9vkjv6njbih9t69t3wfg}.
%

\Subsection{Proof of Theorem~\ref{Llin3maj}}
As previously, assume \( \supH = 1 \).
We use \( \zq(\BBH,\xxc) \leq \sqrt{\dimH} + \sqrt{2 \xxc} \).
Further, \( \constg^{-1} \xxc - \sqrt{2\xxc} + \constg/\sqrt{2} \geq 0 \) and thus,
\begin{EQA}
	\sqrt{2\xxc} - \constg^{-1} \xxc
	& \leq &
	\constg/\sqrt{2} \, .
\label{chf63ghfnb74hrmbgogh}
\end{EQA}
Therefore, for \( \xx \geq \xxc \), it holds 
\begin{EQA}
	\zqc(\BBH,\xx)
	& = &
	\zq(\BBH,\xxc) + \frac{\xx - \xxc}{\constg}
	\leq 
	\sqrt{\dimH} + \sqrt{2 \xxc} - \frac{\xxc}{\constg} + \frac{\xx}{\constg}
	\leq 
	\sqrt{\dimH} + \frac{\constg}{\sqrt{2}} + \frac{\xx}{\constg} \, .
\label{hvtwhnfkdhjwifmbye}
\end{EQA}
In the zone \( \xx \leq \xxc \), it holds \( \zqc(\BBH,\xx) = \zq(\BBH,\xx)	\leq \sqrt{\dimH} + \sqrt{2 \xx} \)
and it remains to note that
\( \sqrt{2 \xx} \leq \constg/\sqrt{2} + \constg^{-1} \xx \).

\Subsection{Proof of Theorem~\ref{TQPlargex}}
Assume w.o.l.g. \( \supH = 1 \).
First consider \( \zq \geq \zqc \).
By \eqref{9fj36vhje37gjgr33erfg} of Theorem~\ref{TQPxivlarge}, 
it holds with \( \constg = \gmb \sqrt{\rhoH} /(2 + \sqrt{\rhoH}) \) 
and \( \xxc = (\zqc - \sqrt{\dimH})^{2}/2 \)
\begin{EQA}
	\P\bigl( \| \QP \xiv \| \geq \zq \bigr)
	& = &
	\P\bigl( \| \QP \xiv \| \geq \zqc + \zq - \zqc \bigr)
	\leq 
	3 \ex^{- \xxc - \constg (\zq - \zqc)} \, .
\label{jcu7wj3fibh99y6ky6hbgi8ed}
\end{EQA}
In particular, \( \P( \| \QP \xiv \| \geq \zqc) \leq 3 \ex^{-\xxc} \).
Integration by parts yields for \( \nuH < \constg \)
\begin{EQA}
	&& \nquad
	\E \ex^{ \nuH (\| \QP \xiv \| - \zqc) } \Ind(\| \QP \xiv \| > \zqc)
	= 
	- \int_{\zqc}^{\infty} \ex^{\nuH (\zq - \zqc)} d\P(\| \QP \xiv \| \geq \zq)
	\\
	&=&
	\P(\| \QP \xiv \| \geq \zqc)
	+ \nuH \int_{\zqc}^{\infty} \ex^{\nuH (\zq - \zqc)} \, \P(\| \QP \xiv \| \geq \zq) \, d\zq 
	\\
	& \leq &
	3 \ex^{-\xxc} + \nuH \int_{\zqc}^{\infty} \ex^{\nuH (\zq - \zqc)- \xxc - \constg (\zq - \zqc)} \, d\zq
	=
	\left( 3 + \frac{3 \nuH}{\constg - \nuH} \right) \, \ex^{-\xxc} 
	\qquad
\label{hgdft6gh3wef7u7ruj543}
\end{EQA}
and \eqref{jhf7yehruybyrhe3wevire} follows.
Similarly, for \( \zq \geq \zqc \), we derive  \eqref{jhf7yehruybyrhe3wevire2} as follows
\begin{EQA}
	&& \nquad
	\E \ex^{ \nuH \| \QP \xiv \| } \Ind(\| \QP \xiv \| > \zq)
	= 
	- \int_{\zq}^{\infty} \ex^{\nuH t} d\P(\| \QP \xiv \| \geq t)
	\\
	& \leq &
	3 \ex^{\nuH \zqc - \xxc - \constg (\zq - \zqc)} 
	+ \frac{3 \nuH}{\constg - \nuH} \, \ex^{\nuH \zqc - \xxc - \constg (\zq - \zqc) } 
	=
	\frac{3 \constg}{\constg - \nuH} \, \ex^{\nuH \zqc -\xxc - (\constg - \nuH) (\zq - \zqc)} \, .
\label{hgdft6gmmm7u7ruj543}
\end{EQA}
Now fix \( \zqx \) with \( \zqx - \sqrt{\dimH} \geq 2 \nuH \) but \( \zqx \leq \zqc \).
Then 
\begin{EQA}
	&& \nquad
	\E \ex^{\nuH \| \QP \xiv \|} \Ind( \| \QP \xiv \| > \zqx)
	=
	- \int_{\zqx}^{\infty} \ex^{\nuH \zq} d\P(\| \QP \xiv \| \geq \zq)
	\\
	&=&
	\ex^{\nuH \zqx} \P(\| \QP \xiv \| \geq \zqx) 
	+ \nuH \left( \int_{\zqx}^{\zqc} + \int_{\zqc}^{\infty} \right) \ex^{\nuH \zq } \P(\| \QP \xiv \| \geq \zq) d\zq \, .
\label{jcy603fgy7e34e7yf8reyhu7ewg}
\end{EQA}
By \eqref{PxivbzzBBroB}, for any \( \zq \in [\zqx , \zqc] \),
it holds in view of \( \zq(\BBH,\xx) \leq \sqrt{\dimH} + \sqrt{2 \xx} \)
\begin{EQA}
	\P(\| \QP \xiv \| \geq \zq)
	& \leq &
	3 \ex^{- (\zq - \sqrt{\dimH})^{2}/2} .
\label{fduyr2ytvfiuh23}
\end{EQA}
As \( (\nuH \zq - (\zq - \sqrt{\dimH})^{2}/2)' = \nuH - \zq + \sqrt{\dimH} \leq - \nuH \) 
for \( \zq - \sqrt{\dimH} \geq 2 \nuH \), it holds
\begin{EQA}
	\nuH \int_{\zqx}^{\zqc} \ex^{\nuH \zq - (\zq - \sqrt{\dimH})^{2}/2 } d\zq
	& \leq &
	\ex^{ \nuH \zqx - (\zqx - \sqrt{\dimH})^{2}/2 } \, \nuH \int_{\zqx}^{\zqc} \ex^{- \nuH (\zq - \zqx) } d\zq 
	\leq 
	\ex^{ \nuH \zqx - (\zqx - \sqrt{\dimH})^{2}/2 } \, 
\label{gtw8ibh9t659o2wmljw8}
\end{EQA}
and also \( \nuH \zqx - (\zqx - \sqrt{\dimH})^{2}/2 > \nuH \zqc - (\zqc - \sqrt{\dimH})^{2}/2 \).
Putting this together with the above bound on \( \int_{\zqc}^{\infty} \ex^{\nuH \zq } \P(\| \QP \xiv \| \geq \zq) d\zq  \) 
as in \eqref{hgdft6gh3wef7u7ruj543}
completes the proof of \eqref{ufikwk3ei9vgkj4k4giw3wl}.


\Section{Deviation bounds for Bernoulli vector sums}
\label{SdevBern}
Let \( Y_{i} \) be independent \( \Bernoulli(\thetas_{i}) \), \( i = 1,\ldots,\nsize \).
We denote \( \Yv = (Y_{i}) \in \R^{\nsize} \).
Weighted sums of the \( Y_{i} \) naturally appear in various statistical tasks 
including classification, binary response models, logistic regression etc. 
Recent applications include e.g. stochastic block modeling; see e.g. \cite{GaZh2017}, \cite{Abbe2018} and references therein,
or ranking from pairwise comparison \cite{ChGa2022} among many others.
We show how the general bounds of Section~\ref{Sdevboundexp} can be used for vector sums of Bernoulli r.v.s.
For a linear mapping \( \Psiv \colon \R^{\nsize} \to \R^{\dimp} \), define
\( \xiv = \Psiv (\Yv - \E \Yv) \).
Below we state some deviation bounds on the squared norm \( \| \xiv \|^{2} \) starting from the univariate case. 

\Subsection{Weighted sums of Bernoulli r.v.'s: univariate case}
Given a collections of weights \( (\weight_{i}) \), define 
\begin{EQA}
	S
	&=&
	\sum_{i=1}^{\nsize} Y_{i} \weight_{i} \, ,
	\\
	\VP^{2}
	&=&
	\Var(S)
	=
	\sum_{i=1}^{\nsize} \thetas_{i}(1-\thetas_{i}) \weight_{i}^{2} \, ,
	\\
	\weights 
	&=& 
	\max_{i} |\weight_{i}| .
\label{Het2Fetetp}
\end{EQA}
First we state a deviation bound for a centered sum \( S - \E S \).

\begin{proposition}
\label{LBeBvM}
Let \( Y_{i} \) be independent \( \Bernoulli(\thetas_{i}) \)
and \( \weight_{i} \in \R \), \( i = 1,\ldots,\nsize \). 
Then \( S = \sum_{i=1}^{\nsize} Y_{i} \weight_{i} \) satisfies
\begin{EQA}
	\log \E \exp\Bigl\{ \frac{\lambda (S - \E S)}{\VP} \Bigr\}
	& \leq &
	\lambda^{2}	,
	\qquad
	\lambda \leq \frac{\log(2) \VP}{\weights} \, .
\label{llEelY1tsel}
\end{EQA}
Furthermore, suppose that given \( \xx \geq 0 \), 
\begin{EQA}
	\VP 
	& \geq &
	\frac{3}{2} \, \weights \sqrt{\xx} \, .
\label{Het2Fetetpzz}
\end{EQA}
Then
\begin{EQA}
	\P\bigl( \VP^{-1} |S - \E S| \geq 2 \sqrt{\xx} \bigr)
	& \leq &
	2 \ex^{-\xx} .
\label{PHetm1SESs2x}
\end{EQA}
Without \eqref{Het2Fetetpzz}, the bound \eqref{PHetm1SESs2x} applies
with \( \VP \) replaced by \( \VP_{\xx} = \VP \vee (3 \, \weights \sqrt{\xx} \, /2) \).
\end{proposition}

\begin{proof}
Without loss of generality assume \( \weights = 1 \), otherwise just rescale all the weights by the factor \( 1/\weights \).
We use that 
\begin{EQA}
	f(u)
	\eqdef
	\log \E \exp\Bigl\{ u (S - \E S) \Bigr\}
	&=&
	\sum_{i=1}^{\nbin} \Bigl[ \log\bigl( \thetas_{i} \ex^{u \weight_{i}} + 1 - \thetas_{i} \bigr) 
	- u \weight_{i} \thetas_{i} \Bigr] \, .
\label{llEelY1tsela}
\end{EQA}
This is an analytic function of \( u \) for \( |u| \leq \log 2 \) satisfying 
\( f(0) = 0 \), \( f'(0) = 0 \), and, with \( \upss_{i} = \log \thetas_{i} - \log(1-\thetas_{i}) \),
\begin{EQA}
	f''(u)
	&=&
	\sum_{i=1}^{\nbin} 
	\frac{\weight_{i}^{2} \, \thetas_{i}(1-\thetas_{i}) \, \ex^{u \weight_{i}}}
		 {( \thetas_{i} \ex^{u \weight_{i}} + 1 - \thetas_{i})^{2}}
	=
	\sum_{i=1}^{\nbin} \frac{\weight_{i}^{2} \, \ex^{\upss_{i}+u \weight_{i}}}{(\ex^{\upss_{i}+u \weight_{i}} + 1)^{2}} 
	=
	\sum_{i=1}^{\nbin} \theta_{i}(u) \bigl\{ 1-\theta_{i}(u) \bigr\} \weight_{i}^{2}
\label{nu2He2fpput}
\end{EQA}
for \( \theta_{i}(u) = \ex^{\upss_{i}+u \weight_{i}}/(\ex^{\upss_{i}+u \weight_{i}} + 1) \).
Clearly \( \theta_{i}(u) \) and thus, \( \theta_{i}(u) \bigl\{ 1-\theta_{i}(u) \bigr\} \) monotonously increases with \( u \) and 
it holds for \( \thetas_{i} = \theta_{i}(0) \) 
\begin{EQA}
	\theta_{i}(u) \bigl\{ 1-\theta_{i}(u) \bigr\}
	& \leq &
	\ex^{|u|} \, \thetas_{i} \, (1 - \thetas_{i})
	\leq 
	2 \, \thetas_{i} \, (1 - \thetas_{i}) ,
	\qquad
	|u| \leq \log 2.
\label{2ts1tseu1tb}
\end{EQA}
This yields 
\begin{EQA}
	f(u)
	& \leq & 
	\VP^{2} \, u^{2}
	\qquad
	|u| \leq \log 2.
\label{flaHetnu2la22}
\end{EQA}
As \( \xx \leq 4\VP^{2}/9 \), the value \( \lambda = \sqrt{\xx} \) fulfills 
\( \lambda/\VP = \sqrt{\xx}/\VP \leq  \log 2 \leq 2^{-1/2} \).
Now by the exponential Chebyshev inequality 
\begin{EQA}
	\P\Bigl( \VP^{-1}(S - \E S) \geq 2 \sqrt{\xx} \Bigr)
	& \leq &
	\exp\bigl\{ - 2 \lambda \sqrt{\xx} + f(\lambda/\VP) \bigr\}
	\\
	& \leq &
	\exp\bigl( - 2 \lambda \sqrt{\xx} + \lambda^{2} \bigr)
	=
	\ex^{-\xx} .
\label{exmxxsq2xnu2la22}
\end{EQA}
Similarly one can bound \( \E S - S \).
\end{proof}

\Subsection{Deviation bounds for Bernoulli vector sums}

Now we present an upper bound on the norm of a vector \( \xiv = \Psiv (\Yv - \E \Yv) \),
where \( \Psiv \) is a linear mapping \( \Psiv \colon \R^{\nsize} \to \R^{\dimp} \).
It holds 
\begin{EQA}
	\Var(\xiv)
	&=&
	\Var(\Psiv \Yv)
	=
	\Psiv \Var(\Yv) \Psiv^{\T} .
\label{HP2VxiWWYT}
\end{EQA}
We aim at bounding the squared norm \( \| \QP \xiv \|^{2} \) for another linear mapping \( \QP \colon \R^{\dimp} \to \R^{\dimq} \).

\begin{theorem}
\label{LHDxibound}
Let \( Y_{i} \sim \Bernoulli(\thetas_{i}) \), \( i = 1,\ldots, n \).
Consider \( \xiv = \Psiv (\Yv - \E \Yv) \), and let \( \HVB^{2} \geq 2 \Var(\xiv) \).
Define 
\begin{EQA}
	\weights
	&=&
	\max_{i \leq \nsize} \| \HVB^{-1} \Psiv_{i} \| \, ,\qquad
	\gmb
	=
	\log(2) / \weights \, .
\label{WHm1uinfmiwi}
\end{EQA} 
Then with \( \BBH = \QP \HVB^{2} \QP^{\T} \) and \( \zqc(\BBH,\xx) \) from \eqref{dyv6ejf8gjwkerih83}, it holds
\begin{EQA}
    \P\bigl( \| \QP \xiv \| \ge \zqc(\BBH,\xx) \bigr)
    & \le &
    3 \ex^{-\xx} .
\label{PxivbzzBBroB3m}
\end{EQA}    
\end{theorem}

\begin{proof}
We apply the general result of Corollary~\ref{CTQPxivlarge} under conditions \eqref{expgamgm}.
For any vector \( \uv \), consider the scalar product 
\( \langle \HVB^{-1} \xiv,\uv \rangle = \langle \HVB^{-1} \Psiv (\Yv - \E \Yv), \uv \rangle \).
It is obviously a weighted centered sum of the Bernoulli r.v.'s \( Y_{i} - \thetas_{i} \) with
\begin{EQA}
	\Var \langle \HVB^{-1} \xiv,\uv \rangle
	& \leq &
	\| \uv \|^{2}/2 .
\label{VlHm1xiuu22}
\end{EQA}
One can write with \( \eps_{i} = Y_{i} - \thetas_{i} \) and \( \epsv = (\eps_{i}) \)
\begin{EQA}
	\langle \HVB^{-1} \xiv,\uv \rangle
	&=&
	\bigl\langle \epsv, \Psiv^{\T} \HVB^{-1} \uv \bigr\rangle .
\label{lHm1xiuepWH1u}
\end{EQA}
By the Cauchy-Schwarz inequality, it holds 
\begin{EQA}
	\| \Psiv^{\T} \HVB^{-1} \uv \|_{\infty}
	=
	\max_{i} \bigl| (\HVB^{-1} \Psiv_{i})^{\T} \uv \bigr|
	& \leq &
	\weights \| \uv \| .
\label{WHm1uinfmiwi}
\end{EQA} 
Bound \eqref{llEelY1tsel} of Proposition~\ref{LBeBvM} on the exponential moments of \( \langle \HVB^{-1} \xiv,\uv \rangle \) implies 
\begin{EQA}
	\log \E \exp\bigl\{ \langle \HVB^{-1} \xiv,\uv \rangle \bigr\}
	& \leq &
	\| \uv \|^{2}/2	,
	\qquad
	\| \uv \| \leq \log(2) / \weights \, .
\label{llEelY1tsel2}
\end{EQA}
Therefore, \eqref{expgamgm} is fulfilled with 
\( \gmb = \log(2) / \weights \). 
The deviation bound \eqref{uyfuyerd7e7uhh8yy689t} of Corollary~\ref{CTQPxivlarge} yields the assertion.
\end{proof}


\def\alpHm{\rho}
\def\Gaussb{\bar{\Gaussv}}
\def\Gauss{\zeta}
\def\Gaussv{\bb{\Gauss}}
\def\BHT{\mathcal{B}}
\def\BBHt{\tilde{\BBH}}
\def\Sigmat{\tilde{\Sigma}}

\Section{Frobenius norm losses for empirical covariance}
\label{SFrobnorm}

Let \( \Xv_{i} \sim \ND(0,\Sigma) \) be i.i.d. zero mean Gaussian vectors in \( \R^{\dimp} \) 
with a covariance matrix \( \Sigma \in \Matr_{\dimp} \).
By \( \Sigmah \) we denote the empirical covariance
\begin{EQA}
	\Sigmah
	& \eqdef &
	\frac{1}{n} \sumi \Xv_{i} \, \Xv_{i}^{\T} \, .
\label{h7hij476wtdysuwjfjti}
\end{EQA}
Our goal is to establish sharp dimension free deviation bounds on the squared Frobenius norm 
\( \| \Sigmah - \Sigma \|_{\Fr}^{2} \):
\begin{EQA}
	\| \Sigmah - \Sigma \|_{\Fr}^{2}
	&=&
	\tr (\Sigmah - \Sigma)^{2}  .
\label{trqtyc6q67ed737duyq}
\end{EQA}
The well developed random matrix theory mainly focuses on the spectral or operator norm of 
\( \Sigmah - \Sigma \); see e.g. \cite{Tropp2015}, \cite{Vers2018} and references therein.
The Frobenius loss \eqref{trqtyc6q67ed737duyq} is much less studied. 
We mention \cite{Lou2014} for the Gaussian case and \cite{BuXi2015} for \( \Xv_{i} \) sub-gaussian.
In some statistical problem like high dimensional random design regression \cite{ChMo2022}, \cite{baLoLu2020} 
or error-in-operator models \cite{Sp2023a}
Frobenius norm of \( \Sigmah - \Sigma \) arises in a natural way.
\cite{ChMo2022} applies Hanson-Wright approach.
We demonstrate how the general results of Section~\ref{Sdevboundexp} can be used for obtaining accurate deviation bounds 
for \( \| \Sigmah - \Sigma \|_{\Fr}^{2} \) and for supporting the concentration phenomenon.

\Subsection{Upper bounds}
First we establish a tight upper bound on \( \| \Sigmah - \Sigma \|_{\Fr}^{2} \).
We identify the matrix \( \Sigmah \) with the vector in the linear subspace of \( \R^{\dimp \times \dimp} \) 
composed by symmetric matrices.
Our aim is in showing that the quantiles of \( \| \Sigmah - \Sigma \|_{\Fr}^{2} \) mimic well similar quantiles 
of \( \| \Sigmat - \Sigma \|_{\Fr}^{2} \) for a Gaussian matrix \( \Sigmat \) with the same covariance structure as \( \Sigmah \).
Define 
\begin{EQA}
	\dimH(\Sigma) 
	= (\tr \Sigma)^{2} + \tr \Sigma^{2} ,
	\qquad
	\vH^{2}(\Sigma)
	&=&
	\bigl( \tr \Sigma^{2} \bigr)^{2} + \tr \Sigma^{4} .
\label{9ckd63hggy7r67y3ghft63}
\end{EQA}
Later we show that \( \dimH(\Sigma) = \E \| \Sigmah - \Sigma \|_{\Fr}^{2} = \tr \Var(\Sigmat) \) and 
\( \vH^{2}(\Sigma) = \tr \{ \Var(\Sigmat) \}^{2} \) while \( \supH(\Sigma) = \| \Var(\Sigmat) \| = 2 \| \Sigma \|^{2} \).
In our results we implicitly assume a high dimensional situation with \( \dimH(\Sigma) \) large.
The presented bounds also require that \( n \gg \dimH(\Sigma) \).

\begin{theorem}
\label{TFrobGauss}
Assume \( \| \Sigma \| = 1 \) and \( \dimH(\Sigma) < n/8 \).
Given \( \xx \) with \( 4 \sqrt{\xx} < \sqrt{n/8} - \sqrt{\dimH(\Sigma)} \),
fix \( \alpHm < 1 \) by 
\begin{EQA}
	\alpHm (1 - \alpHm) \sqrt{n/8} 
	&=&
	\sqrt{\dimH(\Sigma)} + 4 \sqrt{\xx} \, .
\label{yf7j3e53thnutjkerkin}
\end{EQA}
Then 
\begin{EQA}
	\P\Bigl( 
		n \| \Sigmah - \Sigma \|_{\Fr}^{2} > \frac{1}{1 - \alpHm} \bigl\{ \dimH(\Sigma) + 2 \vH(\Sigma) \sqrt{\xx} + 4 \xx \bigr\} 
	\Bigr)
	& \leq &
	3 \ex^{-\xx} \, .
\label{uc6d5e4w4ew5dtdbiie}
\end{EQA}
\end{theorem}


\Subsection{Lower bounds}
This section presents a lower bound on the Frobenius norm of \( \Sigmah - \Sigma \).
Later in Section~\ref{SFrobconc} we state the concentration phenomenon for \( \| \Sigmah - \Sigma \|_{\Fr}^{2} \).

\begin{theorem}
\label{TFrobGausslo}
Let \( \| \Sigma \| = 1 \) and  \( \dimH(\Sigma) \) and \( \vH(\Sigma) \) be defined by \eqref{9ckd63hggy7r67y3ghft63}.
For \( \xx > 0 \) with \( 2 \sqrt{\xx} \leq \dimH(\Sigma)/\vH(\Sigma) \), define \( \muH = \muH(\xx) = 2 \sqrt{\xx}/\vH(\Sigma) \) and assume that
there is \( \alpH < 1/2 \) satisfying
\begin{EQA}
	\alpH \sqrt{\frac{1 - 2\alpH}{1 - \alpH}}
	& \geq &
	\sqrt{\frac{\muH(\xx)}{n}} \biggl( \sqrt{2 \dimH(\Sigma)} + \frac{\sqrt{2} \, \dimH(\Sigma)}{\vH(\Sigma)} \biggr) .
\label{ujw3jdfcv7823ujfwqyshf}
\end{EQA}
Then
\begin{EQA}
	\P\Bigl( n \| \Sigmah - \Sigma \|_{\Fr}^{2} 
		< \frac{1 - 2 \alpH}{1 - \alpH} \, \dimH(\Sigma) - 2 \vH(\Sigma) \sqrt{\xx} 
	\Bigr)
	& \leq &
	2 \ex^{-\xx} .
\label{9fhwe3hdfv7ye43kbhnikj}
\end{EQA}
\end{theorem}

\Subsection{Concentration of the Frobenius loss}
\label{SFrobconc}
Putting together Theorem~\ref{TFrobGauss} and Theorem~\ref{TFrobGausslo} yields the following corollary.

\begin{corollary}
\label{CTFrobGauss}
Under conditions of Theorem~\ref{TFrobGauss} and Theorem~\ref{TFrobGausslo}, 
it holds for any \( \xx \) resolving \eqref{yf7j3e53thnutjkerkin} and \eqref{ujw3jdfcv7823ujfwqyshf}
on a random set \( \Omega(\xx) \) with \( \P\bigl( \Omega(\xx) \bigr) \geq 1 - 5 \ex^{-\xx} \)
\begin{EQA}
	\frac{1 - 2\alpH}{1 - \alpH} \, \dimH(\Sigma) - 2 \vH(\Sigma) \sqrt{\xx} 
	& \leq &
	n \| \Sigmah - \Sigma \|_{\Fr}^{2} 
	\leq  
	\frac{1}{1 - \alpHm} \bigl\{ \dimH(\Sigma) + 2 \vH(\Sigma) \sqrt{\xx} + 4 \xx \bigr\} .
	\qquad
\label{uc6d5e4w4ew5dtdbiiec}
\end{EQA}
\end{corollary}

This result mimics similar bound of Theorem~\ref{TexpbLGA} for \( \Sigmah \) Gaussian 
and of Theorem~\ref{Tdevboundgm} for \( \Sigmah \) sub-Gaussian.
However, the empirical covariance \( \Sigmah \) is quadratic in the \( \Xv_{i} \)'s and thus, only sub-exponential.
We pay an additional factor \( (1 - \alpHm)^{-1} \) in the upper quantile function 
and the factor \( \frac{1 - 2\alpH}{1 - \alpH} \) in the lower quantile function for this extension.

Further we discuss the concentration phenomenon for the Frobenius error \( n \| \Sigmah - \Sigma \|_{\Fr}^{2} \) 
around its expectation \( \dimH(\Sigma) \).
Even in the Gaussian case, it meets only in high-dimensional situation with \( \dimH(\Sigma) \) large.
As \( \vH^{2}(\Sigma) \leq \dimH(\Sigma) \supH(\Sigma) = 2 \dimH(\Sigma) \), this also implies
\( \vH(\Sigma) \ll \dimH(\Sigma) \).
Statement \eqref{uc6d5e4w4ew5dtdbiiec} can be rewritten as
\begin{EQA}
	- \frac{\alpH \, \dimH(\Sigma)}{1 - \alpH} - 2 \vH(\Sigma) \sqrt{\xx} 
	& \leq &
	n \| \Sigmah - \Sigma \|_{\Fr}^{2} - \dimH(\Sigma)
	\leq  
	\frac{\alpHm \, \dimH(\Sigma)}{1 - \alpHm} + \frac{2 \vH(\Sigma) \sqrt{\xx} + 4 \xx}{1 - \alpHm} \, .
	\qquad
\label{uc6d5e4w4ew5dtdbiiecm}
\end{EQA}
Therefore, concentration effect of the loss \( n \| \Sigmah - \Sigma \|_{\Fr}^{2} \) requires \( \dimH(\Sigma) \) large and
\( \alpH \) and \( \alpHm \) small.
Then for \( \xx \ll \dimH(\Sigma) \), quantiles of \( n \| \Sigmah - \Sigma \|_{\Fr}^{2} - \dimH(\Sigma) \)
are smaller in order than \( \dimH(\Sigma) \).
Definition \eqref{yf7j3e53thnutjkerkin} of \( \alpHm \) ensures \( \alpHm \asymp \sqrt{\dimH(\Sigma)/n} \),
and hence, ``\( \alpHm \ll 1 \)'' is equivalent to ``\( \dimH(\Sigma) \ll n \)''. 
Condition ensuring \( \alpH \ll 1 \) is similar.
To see this, assume  \( \vH^{2}(\Sigma) \asymp \dimH(\Sigma) \).
Then \( \xx \ll \dimH(\Sigma) \) yields \( \muH(\xx) = 2 \sqrt{\xx}/\vH(\Sigma) \ll 1 \)
and definition \eqref{ujw3jdfcv7823ujfwqyshf} of \( \alpH \) implies 
\begin{EQA}
	\alpH
	& \lesssim &
	\sqrt{\frac{\muH}{n}} \biggl( \sqrt{2 \dimH(\Sigma)} + \frac{\sqrt{2} \, \dimH(\Sigma)}{\vH(\Sigma)} \biggr)
	\lesssim 
	\sqrt{\frac{\dimH(\Sigma)}{n}} \, .
\label{ojju8733erfsy7r745w23e3r}
\end{EQA}

\Subsection{Weighted Frobenius norm}
The result can be easily extended to the case of a weighted Frobenius norm.
Consider for any linear mapping \( \KH \colon \R^{\dimp} \to \R^{\dimq} \) the value 
\( n \| \KH (\Sigmah - \Sigma) \KH^{\T} \|_{\Fr}^{2} \).

\begin{theorem}
\label{TFrobGaussA}
Let \( \| \Sigma \| = 1 \) and \( \KH \colon \R^{\dimp} \to \R^{\dimq} \) be a linear operator with 
\( \| \KH \| = \| \KH^{\T} \KH \| = 1 \).
Define \( \Sigma_{\KH} \eqdef \KH \Sigma \KH^{\T} \),
\begin{EQA}
	\dimH_{\KH} 
	& \eqdef & 
	\dimH(\Sigma_{\KH})
	=
	\tr^{2} (\Sigma_{\KH}) + \tr (\Sigma_{\KH})^{2},
	\qquad
	\vH_{\KH}^{2} 
	\eqdef 
	\vH^{2}(\Sigma_{\KH})
	=
	\bigl\{ \tr (\Sigma_{\KH}^{2}) \bigr\}^{2} + \tr (\Sigma_{\KH})^{4},
\label{0j3wgyrfg76hy5tjhtdeert}
\end{EQA}
and assume \( \dimH_{\KH} < n/8 \).
The the statements of Theorem~\ref{TFrobGauss} and Theorem~\ref{TFrobGausslo} apply to 
\( n \| \KH (\Sigmah - \Sigma) \KH^{\T} \|_{\Fr}^{2} \) after replacing 
\( \dimH(\Sigma) \) and \( \vH(\Sigma) \) with \( \dimH_{\KH} \) and \( \vH_{\KH} \).
\end{theorem}

\begin{proof}
We can represented  
\begin{EQA}
	\sqrt{n} \, \KH (\Sigmah - \Sigma) \KH^{\T}
	&=&
	\KH \, \Sigma^{1/2} \errS \, \Sigma^{1/2} \KH^{\T}
\label{on0njm9yu7i8erwhedc}
\end{EQA}
with \( \errS \) from \eqref{jufc6jbjuhhjiutifg}.
This reduces the result to the previous case with \( \Sigma_{\KH} = \KH \Sigma \KH^{\T} \) in place of \( \Sigma \).
\end{proof}

\Subsection{Proof of Theorem~\ref{TFrobGauss}}
Each vector \( \gaussv_{i} = \Sigma^{-1/2} \Xv_{i} \) is standard normal.
Define 
\begin{EQA}
	\errS
	&=&
	\frac{1}{n^{1/2}} \sumi (\gaussv_{i} \gaussv_{i}^{\T} - \Id_{\dimp}) .
\label{jufc6jbjuhhjiutifg}
\end{EQA}
We will use the representation \( \Sigmah - \Sigma = n^{-1/2} \Sigma^{1/2} \, \errS \, \Sigma^{1/2} \) and
\begin{EQA}
	n \| \Sigmah - \Sigma \|_{\Fr}^{2}
	&=&
	\tr ( \Sigma^{1/2} \errS \, \Sigma \, \errS \, \Sigma^{1/2} )
	=
	\| \Sigma^{1/2} \errS \, \Sigma^{1/2} \|_{\Fr}^{2} \, .
\label{hvcftw3ytet655rdffhkfj}
\end{EQA}
The main step is in applying Theorem~\ref{Tdevboundgm} to the quadratic form \( \| \QP \errS \|_{\Fr}^{2} \) 
with \( \QP \errS = \Sigma^{1/2} \, \errS \, \Sigma^{1/2} \). 
First check \eqref{expgamgm} for \( \xiv = \errS \).

\begin{lemma}
\label{LexpmgamG}
For any symmetric \( \Gamma \in \Matr_{\dimp} \) with \( \| \Gamma \|_{\Fr} \leq \gmn < \sqrt{n}/2 \), it holds
\begin{EQA}
\label{vfrw44d4d57y8hjhqtstik2}
	\E \langle \Gamma, \errS \rangle^{2}
	&=&
	2 \| \Gamma \|_{\Fr}^{2} \, ,
	\\
	\log \E \exp \langle \Gamma, \errS \rangle
	& \leq &
	\frac{1}{1 - 2 n^{-1/2} \| \Gamma \|} \, \| \Gamma \|_{\Fr}^{2}
	\leq 
	\frac{1}{1 - 2 n^{-1/2} \gmn} \, \| \Gamma \|_{\Fr}^{2} \, .
\label{vfrw44d4d57y8hjhqtstik}
\end{EQA}
\end{lemma}

\begin{proof}
Let us fix any symmetric \( \Gamma \in \Matr_{\dimp} \) with \( \| \Gamma \|_{\Fr} \leq \gmn \).
For the scalar product \( \langle \Gamma, \errS \rangle \), we use the representation
\begin{EQA}
	\langle \Gamma, \errS \rangle
	&=&
	\tr (\Gamma \errS)
	=
	\frac{1}{n^{1/2}} \sumi \bigl\{ \gaussv_{i}^{\T} \Gamma \gaussv_{i} - \E (\gaussv_{i}^{\T} \Gamma \gaussv_{i}) \bigr\} .
\label{5twy7fuuf37dhjuubjbms}
\end{EQA}
Then by independence of the \( \gaussv_{i} \)'s and Lemma~\ref{Gaussmoments}, it holds 
\begin{EQA}
	\E \langle \Gamma, \errS \rangle^{2}
	&=&
	\frac{1}{n} \sumi \E \bigl\{ \gaussv_{i}^{\T} \Gamma \gaussv_{i} - \E (\gaussv_{i}^{\T} \Gamma \gaussv_{i}) \bigr\}^{2}
	=
	2 \tr \Gamma^{2} .
\label{5twy7fuuf37dhjuubjbms2}
\end{EQA}
Now consider the exponential moment of \( \langle \Gamma, \errS \rangle \).
Again, independence of the \( \gaussv_{i} \)'s yields
\begin{EQA}
	\log \E \exp \langle \Gamma, \errS \rangle
	&=&
	\sumi \log \E \exp \frac{\gaussv_{i}^{\T} \Gamma \gaussv_{i}}{\sqrt{n}} - \sqrt{n} \, \tr \Gamma
	\\
	&=&
	\frac{n}{2} \log \det \bigl( \Id_{\dimp} - \frac{2}{\sqrt{n}} \, \Gamma \bigr) - \sqrt{n} \, \tr \Gamma
	\qquad
\label{hw6jhydhqanchyenqqdbsk} 
\end{EQA}
provided that \( 2\Gamma < \sqrt{n} \Id_{\dimp} \).
Moreover, by Lemma~\ref{Lqfexpmom}
\begin{EQA}
	\left| \frac{n}{2} \log \det (\Id_{\dimp} - 2 n^{-1/2} \Gamma) - \sqrt{n} \, \tr \Gamma \right|
	& \leq &
	\frac{\tr \Gamma^{2}}{1 - 2 n^{-1/2} \| \Gamma \|} 
	=
	\frac{\| \Gamma \|_{\Fr}^{2}}{1 - 2 n^{-1/2} \| \Gamma \|} \, ,
	\qquad
\label{hyakboj84gjqtqgxjvu}
\end{EQA}
and the assertion follows in view of \( \| \Gamma \| \leq \| \Gamma \|_{\Fr} \leq \gmn \).
\end{proof}


We now fix \( \gmn = \alpHm \sqrt{n}/2 \).
Then the random matrix \( \xiv = \errS \) follows
condition \eqref{expgamgm} with 
\( \HVB^{2} = 2 (1 - \alpHm)^{-1} \Id \).
This enables us to apply Theorem~\ref{Tdevboundgm} to the quadratic form \( \| \QP \errS \|_{\Fr}^{2} \) 
for \( \QP \errS = \Sigma^{1/2} \, \errS \, \Sigma^{1/2} \). 
By \eqref{vfrw44d4d57y8hjhqtstik2}, it holds \( \Var(\errS) = 2 \Id \).
Now introduce a Gaussian element \( \errSt \) with the same covariance structure. 
One can use \( \errSt = (\Gaussv + \Gaussv^{\T})/\sqrt{2} \), where \( \Gaussv = (\Gauss_{ij}) \) is a random \( \dimp \)-matrix 
with i.i.d. standard normal entries \( \Gauss_{ij} \).
Indeed, for any symmetric \( \dimp \)-matrix \( \Gamma \),
\begin{EQA}
	\E \langle \errSt, \Gamma \rangle^{2}
	&=&
	2 \E \langle \Gaussv,\Gamma \rangle^{2}
	=
	2 .
\label{8jjjii8i8334eeee4kfoof}
\end{EQA}
Statement \eqref{PxivbzzBBroB} of Theorem~\ref{Tdevboundgm} yields nearly the same deviation bounds 
for \( \| \QP \errS \|_{\Fr}^{2} \) as for \( \| \QP \errSt \|_{\Fr}^{2} \) with \( \errSt \sim \ND(0,\Var(\errS)) \).
Theorem~\ref{TexpbLGA} claims 
\begin{EQA}
	\P\bigl( \| \QP \errSt \|_{\Fr}^{2} > \zq^{2}(\BBHt,\xx) \bigr)
	& \leq &
	\ex^{-\xx} \, ,
\label{jf67wyfuvb87e7eufhgghwxu}
\end{EQA}
where \( \BBHt = \Var(\QP \errSt) \) and 
the quantile \( \zq(\BBH,\xx) \) is defined as
\begin{EQA}
	\zq^{2}(\BBH,\xx)
	&=&
	\tr \BBH + 2 \sqrt{\xx \tr(\BBH^{2})} + 2 \xx \| \BBH \| .
\label{fiu47fkw3ydfywe3hdf7y2h}
\end{EQA}

\begin{lemma}
\label{LBBHt}
Let \( \errSt = (\Gaussv + \Gaussv^{\T})/\sqrt{2} \), where \( \Gaussv = (\Gauss_{ij}) \) is a random \( \dimp \)-matrix 
with i.i.d. standard normal entries \( \Gauss_{ij} \).
Consider \( \QP \errSt = \Sigma^{1/2} \, \errSt \, \Sigma^{1/2} \).
It holds for \( \BBHt = \Var(\QP \errSt) \)
\begin{EQA}
	\tr \BBHt
	&=&
	\dimH(\Sigma) \, ,
	\quad
	\tr \BBHt^{2}
	=
	\vH^{2}(\Sigma) \, ,
	\quad
	\| \BBHt \|
	=
	2. 
\label{yx7hqwjxc78g94mg6e6e4heg}
\end{EQA}
\end{lemma}

\begin{proof}
We may assume \( \Sigma = \diag\{ \lambda_{1},\ldots,\lambda_{\dimp} \} \). 
Then it holds by Lemma~\ref{Gaussmoments}
\begin{EQA}
	\| \QP \errSt \|_{\Fr}^{2}
	=
	\| \Sigma^{1/2} \, \errSt \, \Sigma^{1/2} \|_{\Fr}^{2}
	=
	\frac{1}{2} \sum_{i,j=1}^{\dimp} \lambda_{i} \, \lambda_{j} \, (\Gauss_{ij} + \Gauss_{ji})^{2}
	& \eqd &
	2 \sum_{i \leq j} \lambda_{i} \, \lambda_{j} \, \Gauss_{ij}^{2}
	\qquad
\label{udhf6e3ggun8true4hfvy}
\end{EQA}
and thus
\begin{EQA}
	\tr \BBHt
	=
	\E \| \QP \errSt \|_{\Fr}^{2}
	&=&
	2 \sum_{i \leq j} \lambda_{i} \, \lambda_{j}
	=
	\biggl( \sum_{i=1}^{\dimp} \lambda_{i} \biggr)^{2} + \sum_{i=1}^{\dimp} \lambda_{i}^{2}
	=
	\dimH(\Sigma) \, .
\label{yx7hqwjxc78wjuanct6whw}
\end{EQA}
Further we compute \( \vH^{2}(\Sigma) = \tr \BBHt^{2} \).
Note that \( \Var( \| \QP \errSt \|_{\Fr}^{2} ) \neq \Var( \| \QP \errS \|_{\Fr}^{2} ) \).
Due to Lemma~\ref{Gaussmoments}, it holds
\( \vH^{2}(\Sigma) = \Var( \| \QP \errSt \|_{\Fr}^{2} )/2 \) yielding by \eqref{udhf6e3ggun8true4hfvy}
\begin{EQA}
	\vH^{2}(\Sigma)
	&=&
	2 \sum_{i \leq j} \lambda_{i}^{2} \lambda_{j}^{2} \Var (\Gauss_{ij}^{2})	
	=
	2 \sum_{i \neq j} \lambda_{i}^{2} \lambda_{j}^{2} 
	+ 2 \sum_{i=1}^{\dimp} \lambda_{i}^{4} 
	=
	\bigl( \tr \Sigma^{2} \bigr)^{2} + \tr \Sigma^{4} .
\label{ioer49u8fdj42e3u7fdjgt}
\end{EQA}
Finally, \( \Var(\errS) = 2 \Id \) and \( \| \Sigma \| = 1 \) implies \( \supH(\Sigma) = \| \QP \Var(\errS) \QP^{\T} \| = 2 \).
\end{proof}

Now we apply Theorem~\ref{Tdevboundgm} to 
\( n \| \Sigmah - \Sigma \|_{\Fr}^{2} = \| \QP \errS \|_{\Fr}^{2} \).
Following to Lemma~\ref{LexpmgamG}, define \( \BBH = (1 - \nuH)^{-1} \BBHt \).
Then with \( \zq^{2}(\BBH,\xx) \) from \eqref{fiu47fkw3ydfywe3hdf7y2h}
\begin{EQA}
	\P\Bigl( n \| \Sigmah - \Sigma \|_{\Fr}^{2} > \zq^{2}(\BBH,\xx) \Bigr)
	& = &
	\P\bigl( \| \QP \errS \|_{\Fr}^{2} > \zq^{2}(\BBH,\xx) \bigr)
	\leq 
	3 \ex^{-\xx} ,
	\quad
	\xx \leq \xxc \, ,
\label{ncye6e3y7fiuh8ytst5q}
\end{EQA}
and assertion \eqref{uc6d5e4w4ew5dtdbiie} follows in view of Lemma~\ref{LBBHt} and 
\( \zq^{2}(\BBH,\xx) = (1 - \nuH)^{-1} \zq^{2}(\BBHt,\xx) \).
However, it is still necessary to check that the upper bound \eqref{uc6d5e4w4ew5dtdbiie} 
applies for a given \( \xx \).
\eqref{if7h3rhgy4676rfhdsjw} provides a sufficient condition 
\( \gmn/\supH \geq \sqrt{\dimH/\supH} + \sqrt{8\xx} \) with \( \dimH = \dimH(\Sigma)/(1 - \alpHm) \) and
\( \supH = 2/(1 - \alpHm) \) for \( \gmn = \alpHm \sqrt{n} / 2 \).
%
By \eqref{yf7j3e53thnutjkerkin}
\begin{EQA}
	\frac{\gmn}{\supH} - \sqrt{\frac{\dimH}{\supH}} 
	& = & 
	\frac{\alpHm \sqrt{n}}{2\supH} - \sqrt{\frac{\dimH(\Sigma)}{2}} 
	\geq 
	\frac{\alpHm (1 - \alpHm) \sqrt{n}}{4} - \sqrt{\frac{\dimH(\Sigma)}{2}}
	> 
	\sqrt{8\xx} 
\label{nctc6yw3f76h7u87y98}
\end{EQA}
and the result follows.

\Subsection{Proof of Theorem~\ref{TFrobGausslo}}
As in the proof of the upper bound, we apply Markov's inequality
\begin{EQA}
	\P\bigl( n \| \Sigmah - \Sigma \|_{\Fr}^{2} < \zq \bigr)
	& \leq &
	\ex^{\muH \zq/2} \E \exp\Bigl( - \frac{\muH}{2} n \| \Sigmah - \Sigma \|_{\Fr}^{2} \Bigr) .
\label{9vjnegf66ruhiti4rhsw000}
\end{EQA}
However, now we are free to choose any positive \( \muH \).
Later we evaluate the exponential moments of \( - n \| \Sigmah - \Sigma \|_{\Fr}^{2} \) for all \( \muH > 0 \)
and then, given \( \xx \), fix \( \muH \) and \( \zq \) similarly to the Gaussian case to ensure 
the prescribed deviation probability \( \ex^{-\xx} \). 

Denote by \( \Gaussv = (\Gauss_{ij}) \) a random \( \dimp \times \dimp \) matrix 
with i.i.d. standard Gaussian entries \( \Gauss_{ij} \) and \( \Gaussb \eqdef (\Gaussv + \Gaussv^{\T})/2 \).
Then for any \( \muH > 0 \)
\begin{EQA}
	\exp \bigl( - \muH n \| \Sigmah - \Sigma \|_{\Fr}^{2}/2 \bigr)
	&=&
	\E_{\Gauss} \exp\bigl\{ \imi \sqrt{\muH n} \, \langle \Sigmah - \Sigma, \Gaussv \rangle \bigr\}
	=
	\E_{\Gauss} \exp\bigl\{ \imi \sqrt{\muH n} \, \langle \Sigmah - \Sigma, \Gaussb \rangle \bigr\} .
\label{9vje6e36tghu4i3bieijd}
\end{EQA}
Therefore, by independence of the \( \Xv_{i} \)'s
\begin{EQA}
\label{yhs7vjeje8bjekwytvya}
	\E \exp \bigl( - \muH n \| \Sigmah - \Sigma \|_{\Fr}^{2}/2 \bigr)
	&=&
	\E_{\Gauss} \E \exp\bigl( \imi \sqrt{\muH n} \, \langle \Sigmah - \Sigma, \Gaussb \rangle \bigr)
	\\
	&=&
	\E_{\Gauss} \Bigl\{ \E \exp\bigl( \imi \sqrt{\muH/n} \, \langle \Xv_{1} \Xv_{1}^{\T} - \Sigma, \Gaussb \rangle \bigr) \Bigr\}^{n}
	\\
	&=&
	\E_{\Gauss} \Bigl\{ \E
	\exp\bigl( 
		\imi \sqrt{\muH/n} \, \langle \gaussv \gaussv^{\T} - \Id_{\dimp} \, , \Sigma^{1/2} \, \Gaussb \, \Sigma^{1/2} \rangle 
	\bigr) 
	\Bigr\}^{n} .
\end{EQA}
Further, by Lemma~\ref{Lqfexpmom}, with \( \BHT = \Sigma^{1/2} \, \Gaussb \, \Sigma^{1/2} \)
\begin{EQA}
	&& \nquad
	\Bigl\{ \E
		\exp\bigl( \imi \sqrt{\muH/n} \, \langle \gaussv \gaussv^{\T} - \Id_{\dimp} \, , \BHT \rangle \bigr) 
	\Bigr\}^{n}
	\\
	&=&
	\exp\bigl\{ n \log \det\bigl( \Id_{\dimp} - 2 \imi \sqrt{\muH/n} \, \BHT \bigr)^{-1/2} 
	- \imi \sqrt{\muH n} \tr (\BHT) \bigr\} .
\label{yhs7vjeje8bjekwytgm}
\end{EQA}
Let some \( \xx > 0 \) and some \( \alpH \in (0,1/2) \) be fixed.
Define
\begin{EQA}
	\muH 
	& \eqdef &
	\frac{2 \sqrt{\xx} }{\vH(\Sigma)} \, ,
	\qquad
	\muHa
	\eqdef
	\frac{1 - \alpH}{1 - 2\alpH} \, \muH
	=
	\frac{1 - \alpH}{1 - 2\alpH} \,\, \frac{2 \sqrt{\xx} }{\vH(\Sigma)} \, ,
\label{nsducf7enhw3e7yfgy6ryer4b}
\end{EQA}
and introduce a random set \( \Omega(\alpH) \) with
\begin{EQA}
	\Omega(\alpH)
	& \eqdef &
	\bigl\{ \Gaussv \colon 2 \sqrt{\muHa/n} \, \| \BHT \|
	\leq 
	\alpH \bigr\} ,
	\qquad
	\BHT = \Sigma^{1/2} \, (\Gaussv + \Gaussv^{\T}) \, \Sigma^{1/2}/2.
\label{uiduuvbinjoko4e7gyuuh2w2}
\end{EQA}
It holds on \( \Omega(\alpH) \) 
by \eqref{yhs7vjeje8bjekwytgm} similarly to \eqref{vbu7j3hg8hryhghyidegwdg} of Lemma~\ref{Lqfexpmom} 
\begin{EQA}
	&& \nquad
	\E^{n} \exp\bigl\{ \imi \sqrt{\muHa/n} \, \langle \gaussv \gaussv^{\T} - \Id_{\dimp}, \BHT \rangle \bigr\}
	\leq 
	\exp \biggl( - \muHa \tr(\BHT^{2}) + \frac{\muHa \, \alpH \tr(\BHT^{2})}{1 - \alpH} \biggr)
	\\
	&=&
	\exp\biggl( - \frac{1 - 2 \alpH}{1 - \alpH} \, \muHa \, \tr(\BHT^{2}) \biggr)
	=
	\exp\bigl( - \muH \, \tr(\BHT^{2}) \bigr) .
	\qquad
\label{yhs7vjeje8behfydrjekwytgm}
\end{EQA}
Exponential moments of \( \tr(\BHT^{2}) \) from \eqref{uc8ewjfbu7fue77egkbire} under \( \P_{\Gauss} \) can be easily computed.
We proceed assuming \( \Sigma = \diag\{ \lambda_{j} \} \) and
using that \( \Gauss_{ij} + \Gauss_{ji} \sim \ND(0,2) \) for \( i \neq j \), and all 
\( \Gauss_{ij} + \Gauss_{ji} \) are mutually independent for \( i \leq j \).
This implies
\begin{EQA}
	\tr(\BHT^{2})
	&=&
	\frac{1}{4} \sum_{i,j=1}^{\dimp} \lambda_{i} \, \lambda_{j} \, (\Gauss_{ij} + \Gauss_{ji})^{2}
	\, \eqd \,
	\sum_{i \leq j} \lambda_{i} \, \lambda_{j} \, \Gauss_{ij}^{2}
\label{uc8ewjfbu7fue77egkbire}
\end{EQA}
and 
\begin{EQA}
	\E_{\Gauss} \tr(\BHT^{2})
	&=&
	\sum_{i \leq j} \lambda_{i} \, \lambda_{j}
	=
	\frac{\dimH(\Sigma)}{2} \, ,
\label{87wkfv8ieigy5e5tgvuf}
	\\
	\E_{\Gauss} \exp\{ - \muH \tr(\BHT^{2}) \}
	&=&
	\E_{\Gauss} \exp\biggl( - \muH \sum_{i \leq j} \lambda_{i} \, \lambda_{j} \, \Gauss_{ij}^{2}  \biggr)
	=
	\exp \biggl( - \frac{1}{2} \sum_{i \leq j} \log (1 + 2 \muH \, \lambda_{i} \, \lambda_{j}) \biggr) .
\label{jcvuyfjhi8u97iew22g}
\end{EQA}
The latter expression can be evaluated by using \eqref{m2v241m41mbn} of Lemma~\ref{Lqfexpmom}:
\begin{EQA}
	\E_{\Gauss} \exp\{ - \muH \tr(\BHT^{2}) \}
	& \leq &
	\exp \biggl( 
		- \muH \sum_{i \leq j} \lambda_{i} \, \lambda_{j} + \muH^{2} \sum_{i \leq j} \lambda_{i}^{2} \, \lambda_{j}^{2} 
	\biggr) 
	=
	\exp\biggl( - \frac{\muH \, \dimH(\Sigma)}{2} + \frac{\muH^{2} \vH^{2}(\Sigma)}{4} \biggr).
\label{jcvuyfjhbir64gtughew22g}
\end{EQA}
This and \eqref{yhs7vjeje8behfydrjekwytgm} yield 
\begin{EQA}
	\E \exp \biggl( - \frac{\muHa}{2} n \| \Sigmah - \Sigma \|_{\Fr}^{2} \biggr)
	& \leq &
	\P_{\Gauss}\bigl( \Omega(\alpH)^{c} \bigr) 
	+ \exp \biggl( - \frac{\muH \, \dimH(\Sigma)}{2} + \frac{\muH^{2} \vH^{2}(\Sigma)}{4} \biggr) 
\label{ifgo9elrgpgjn9u9tkfywfgi}
\end{EQA}
and for any \( \zq \) by Markov's inequality \eqref{9vjnegf66ruhiti4rhsw000}
\begin{EQA}
	\P\bigl( n \| \Sigmah - \Sigma \|_{\Fr}^{2} < \zq \bigr)
	& \leq &
	\ex^{\muHa \zq/2} \P_{\Gauss}\bigl( \Omega(\alpH)^{c} \bigr) 
	+ \exp\biggl( \frac{\muHa \zq}{2} - \frac{\muH \, \dimH(\Sigma)}{2} + \frac{\muH^{2} \vH^{2}(\Sigma)}{4} \biggr) .
\label{9vjnegf66ruhiti4rhsw}
\end{EQA}
With \( \muH = 2 \sqrt{\xx} / \vH(\Sigma) \), we  
define \( \zq \) by
\begin{EQA}
	\muHa \,\zq
	&=&
	\muH \{ \dimH(\Sigma) - 2 \vH(\Sigma) \sqrt{\xx} \} 
	=
	\frac{2 \sqrt{\xx} \, \dimH(\Sigma)}{\vH(\Sigma)} - 4 \xx \, 
\label{8ujek94jhkmokjkoiowndy}
\end{EQA}
yielding 
\begin{EQA}
	\frac{\muHa \zq}{2} - \frac{\muH \, \dimH(\Sigma)}{2} + \frac{\muH^{2} \vH^{2}(\Sigma)}{4}
	&=&
	\frac{\muH}{2} \bigl\{ \dimH(\Sigma) - 2 \vH(\Sigma) \sqrt{\xx} \bigr\}
	- \frac{\muH \, \dimH(\Sigma)}{2} + \frac{\muH^{2} \vH^{2}(\Sigma)}{4}
	=
	-\xx
\label{nyhwnvu7u7w2hvuyghiuhyiu}
\end{EQA}
and
\begin{EQA}
	\P( n \| \Sigmah - \Sigma \|_{\Fr}^{2} < \zq)
	& \leq &
	\ex^{-\xx} + \ex^{\muHa \zq/2} \P_{\Gauss}\bigl( \Omega(\alpH)^{c} \bigr)
\label{y7uemb8jrboboierkedfcoaq}
\end{EQA}
where
\begin{EQA}
	\zq
	&=&
	\biggl( 1 - \frac{\alpH}{1 - \alpH} \biggr) \bigl\{ \dimH(\Sigma) - 2 \vH(\Sigma) \sqrt{\xx} \bigr\}
	\geq 
	\dimH(\Sigma)
	- \frac{\alpH}{1 - \alpH} \, \dimH(\Sigma)
	- 2 \vH(\Sigma) \sqrt{\xx} \, .
\label{kfuhyg6er4y6bfhggfy3}
\end{EQA}
For bounding the probability of the set \( \Omega(\alpH)^{c} \) from \eqref{uiduuvbinjoko4e7gyuuh2w2}, 
one can apply the advanced results from the random matrix theory.
To keep the proof self-contained, we use a simple bound \( \| \BHT \|^{2} \leq \| \BHT \|_{\Fr}^{2} = \tr(\BHT^{2}) \).
For any matrix \( \Gamma \), it holds 
\begin{EQA}
	\Var\langle \Gaussb,\Gamma \rangle 
	&=&
	\frac{1}{4} \E \biggl( \sum_{i,j=1}^{\dimp} \Gamma_{ij} (\Gauss_{ij} + \Gauss_{ji}) \biggr)^{2}
	=
	\| \Gamma \|_{\Fr}^{2} \, 
\label{onhr6g7e8gheh3vujweg}
\end{EQA}
yielding \( \| \Var(\Gaussb) \| \leq 1 \) and \( \| \Var(\BHT) \| \leq 1 \).
Also by \eqref{87wkfv8ieigy5e5tgvuf} \( \E \| \BHT \|_{\Fr}^{2} = \dimH(\Sigma)/2 \).
Therefore, by Theorem~\ref{TexpbLGA} applied to \( \| \BHT \|_{\Fr}^{2} \), it holds for any \( \xxs \)
\begin{EQA}
	\P_{\Gauss} \bigl( \| \BHT \|_{\Fr} > \sqrt{\dimH(\Sigma)/2} + \sqrt{2\xxs} \bigr)
	& \leq &
	\ex^{-\xxs} .
\label{j8jfkg99gue3g7gyw3bjb}
\end{EQA}
By \eqref{nsducf7enhw3e7yfgy6ryer4b} and \eqref{8ujek94jhkmokjkoiowndy}, it holds 
\begin{EQA}
	\xxs
	\eqdef
	\xx + \frac{\muHa \zq}{2} 
	& \leq &
	\frac{\dimH(\Sigma) \sqrt{\xx}}{\vH(\Sigma)} - \xx
	\leq 
	\frac{\dimH^{2}(\Sigma)}{4 \vH^{2}(\Sigma)}
\label{udf73ujerdf8t754rjfcu8wei}
\end{EQA}
and 
\begin{EQA}
	\P_{\Gauss} \biggl( \| \BHT \|_{\Fr} > \sqrt{\frac{\dimH(\Sigma)}{2}} + \frac{\dimH(\Sigma)}{\sqrt{2} \, \vH(\Sigma)} \biggr)
	& \leq &
	\ex^{-\xx - \muHa \zq/2} .
\label{j8jfkg99gvik4r9hijhew3bjb}
\end{EQA}
Therefore, 
by definition \eqref{uiduuvbinjoko4e7gyuuh2w2} and condition \eqref{ujw3jdfcv7823ujfwqyshf}
\begin{EQA}
	\ex^{\muHa \zq/2} \P_{\Gauss}\bigl( \Omega(\alpH)^{c} \bigr)
	& \leq &
	\ex^{\muHa \zq/2} \P\biggl( \| \BHT \|_{\Fr} > \frac{\alpH \sqrt{n}}{2 \sqrt{\muHa}} \biggr)
	\leq 	
	\ex^{-\xx} 
\label{knegyeg67gh47rskerur}
\end{EQA}
and the result follows.

\newpage
\appendix


\Section{Moments of a Gaussian quadratic form}
\label{SmomentqfG}
Let \( \gaussv \) be standard normal in \( \R^{\dimp} \) for \( \dimp \leq \infty \).
Given a self-adjoint trace operator \( \BBH \), consider a quadratic form 
\( \bigl\langle \BBH \gaussv, \gaussv \bigr\rangle \).

\begin{lemma}
\label{Gaussmoments}
It holds
\begin{EQA}
	\E \bigl\langle \BBH \gaussv, \gaussv \bigr\rangle 
	&=& 
	\tr \BBH .
\label{EAarAtrA2trA2}
\end{EQA}
Moreover, 
\begin{EQA}
	\E \bigl( \bigl\langle \BBH \gaussv, \gaussv \bigr\rangle - \tr \BBH \bigr)^{2}
	&=&
	2 \tr \BBH^{2}  ,
	\\
	\E \bigl( \bigl\langle \BBH \gaussv, \gaussv \bigr\rangle - \tr \BBH \bigr)^{3}
	&=&
	8 \tr \BBH^{3} ,
	\\
	\E \bigl( \bigl\langle \BBH \gaussv, \gaussv \bigr\rangle - \tr \BBH \bigr)^{4}
	&=&
	48 \tr \BBH^{4} + 12 (\tr \BBH^{2})^{2} ,
\label{2pG2trD2DGm22m2}
\end{EQA}
and
\begin{EQA}
	\E \bigl\langle \BBH \gaussv, \gaussv \bigr\rangle^{2}
	&=&
	(\tr \BBH)^{2} + 2 \tr \BBH^{2},
	\\
	\E \bigl\langle \BBH \gaussv, \gaussv \bigr\rangle^{3}
	& = &
	(\tr \BBH)^{3} + 6 \tr \BBH \,\, \tr \BBH^{2} + 8 \tr \BBH^{3} ,
	\\
	\E \bigl\langle \BBH \gaussv, \gaussv \bigr\rangle^{4}
	& = &
	(\tr \BBH)^{4} + 12 (\tr \BBH)^{2} \tr \BBH^{2}
	+ 32 (\tr \BBH) \tr \BBH^{3}
	+ 48 \tr \BBH^{4} + 12 (\tr \BBH^{2})^{2} ,
\label{2pG2trD2DGm22m2}
	\\
	\Var \bigl\langle \BBH \gaussv, \gaussv \bigr\rangle^{2}
	& = &
	8 (\tr \BBH)^{2} \tr \BBH^{2}
	+ 32 (\tr \BBH) \tr \BBH^{3}
	+ 48 \tr \BBH^{4} + 8 (\tr \BBH^{2})^{2} .
\label{2pG2trD2DGm22m4}
\end{EQA}
Moreover, if \( \BBH \leq \Id_{\dimp} \) and \( \dimH = \tr \BBH \), then \( \tr \BBH^{m} \leq \dimH \| \BBH \|^{m-1} \) for 
\( m \geq 1 \) and
\begin{EQA}[rcccl]
	\E \bigl\langle \BBH \gaussv, \gaussv \bigr\rangle^{2}
	& \leq &
	\dimH^{2} + 2 \dimH \| \BBH \|
	&\leq &
	(\dimH + \| \BBH \|)^{2},
	\\
	\E \bigl\langle \BBH \gaussv, \gaussv \bigr\rangle^{3}
	& \leq &
	\dimH^{3} + 6 \dimH^{2} \| \BBH \| + 8 \dimH \| \BBH \|^{2}
	&\leq &
	(\dimH + 2 \| \BBH \|)^{3},
	\\
	\E \bigl\langle \BBH \gaussv, \gaussv \bigr\rangle^{4}
	& \leq &
	\dimH^{4} + 12 \dimH^{3} \| \BBH \|
	+ 44 \dimH^{2} \| \BBH \|^{2}
	+ 48 \dimH \| \BBH \|^{3}
	&\leq &
	(\dimH + 3 \| \BBH \|)^{4},
\label{2pG2trD2DGm22m2}
	\\
	\Var \bigl\langle \BBH \gaussv, \gaussv \bigr\rangle^{2}
	& \leq &
	8 \dimH^{3} + 40 \dimH^{2} \| \BBH \| + 48 \dimH \| \BBH \|^{2}.
\label{2pG2trD2DGm22m4}
\end{EQA}
Finally,
\begin{EQA}
	\E (\gaussv \gaussv^{\T} - \Id_{\dimp}) \BBH (\gaussv \gaussv^{\T} - \Id_{\dimp}) 
	&=&
	\BBH + \tr (\BBH) \Id_{\dimp}
\label{njt66777888723fdgy}
\end{EQA}
yielding
\begin{EQA}
	\E \| \BBH (\gaussv \gaussv^{\T} - \Id_{\dimp}) \|_{\Fr}^{2}
	&=&
	(\tr \BBH)^{2} + \tr \BBH^{2} .
\label{njt66777888723fdgyf}
\end{EQA}
\end{lemma}

\begin{proof}
Let \( \chi = \gauss^{2} - 1 \) for \( \gauss \) standard normal.
Then \( \E \chi = 0 \), \( \E \chi^{2} = 2 \), \( \E \chi^{3} = 8 \), \( \E \chi^{4} = 60 \).
Without loss of generality assume \( \BBH \) diagonal: \( \BBH = \diag(\supH_{1},\supH_{2},\ldots,\supH_{\dimp}) \).
Then 
\begin{EQA}
	\xi
	\eqdef
	\bigl\langle \BBH \gaussv, \gaussv \bigr\rangle - \tr \BBH
	&=&
	\sum_{j=1}^{\dimp} \supH_{j} (\gauss_{j}^{2} - 1) ,
\label{j1ljgj2m1}
\end{EQA}
where \( \gauss_{j} \) are i.i.d. standard normal. 
This easily yields
\begin{EQA}
	\E \xi^{2}
	&=&
	\sum_{j=1}^{\dimp} \supH_{j}^{2} \E (\gauss_{j}^{2} - 1)^{2}
	=
	\E \chi^{2} \, \tr \BBH^{2} 
	=
	2 \tr \BBH^{2}  ,
	\\
	\E \xi^{3}
	&=&
	\sum_{j=1}^{\dimp} \supH_{j}^{3} \E (\gauss_{j}^{2} - 1)^{3}
	=
	\E \chi^{3} \, \tr \BBH^{3} 
	=
	8 \tr \BBH^{3} ,
	\\
	\E \xi^{4}
	&=&
	\sum_{j=1}^{\dimp} \supH_{j}^{4} (\gauss_{j}^{2} - 1)^{4}
	+ \sum_{i\neq j} \supH_{i}^{2} \supH_{j}^{2} \E (\gauss_{i}^{2} - 1)^{2} \E (\gauss_{j}^{2} - 1)^{2}
	\\
	&=&
	\bigl( \E \chi^{4} - 3 (\E \chi^{2})^{2} \bigr) \tr \BBH^{4} + 3 (\E \chi^{2} \, \tr \BBH^{2})^{2}
	=
	48 \tr \BBH^{4} + 12 (\tr \BBH^{2})^{2} ,
\label{2pG2trD2DGm22m2}
\end{EQA}
ensuring
\begin{EQA}
	\E \bigl\langle \BBH \gaussv, \gaussv \bigr\rangle^{2}
	&=&
	\bigl( \E \bigl\langle \BBH \gaussv, \gaussv \bigr\rangle \bigr)^{2} 
	+ \E \xi^{2}
	= 
	(\tr \BBH)^{2} + 2 \tr \BBH^{2},
	\\
	\E \bigl\langle \BBH \gaussv, \gaussv \bigr\rangle^{3}
	& = &
	\E \bigl( \xi + \tr \BBH \bigr)^{3}
	=
	(\tr \BBH)^{3} + \E \xi^{3}
	+ 3 \tr \BBH \,\, \E \xi^{2}
	\\
	&=&
	(\tr \BBH)^{3} + 6 \tr \BBH \,\, \tr \BBH^{2} + 8 \tr \BBH^{3} ,
\label{2pG2trD2DGm22m2}
\end{EQA}
and 
\begin{EQA}
	\Var \bigl\langle \BBH \gaussv, \gaussv \bigr\rangle^{2}
	& = &
	\E \bigl( \xi + \tr \BBH \bigr)^{4}
	- \bigl( \E \bigl\langle \BBH \gaussv, \gaussv \bigr\rangle \bigr)^{2}
	\\
	&=&
	\bigl( \tr \BBH \bigr)^{4} + 6 (\tr \BBH)^{2} \E \xi^{2} + 4 \tr \BBH \, \E \xi^{3} + \E \xi^{4}
	- \bigl( (\tr \BBH)^{2} + 2 \tr \BBH^{2} \bigr)^{2}
	\\
	&=& 
	8 (\tr \BBH)^{2} \tr \BBH^{2}
	+ 32 (\tr \BBH) \tr \BBH^{3}
	+ 48 \tr \BBH^{4} + 8 (\tr \BBH^{2})^{2} .
\label{2pG2trD2DGm22m4}
\end{EQA}
For the last result of the lemma, observe that with \( \BBH = \diag(\supH_{1},\supH_{2},\ldots,\supH_{\dimp}) \), 
the matrix \( \E (\gaussv \gaussv^{\T} - \Id_{\dimp}) \BBH (\gaussv \gaussv^{\T} - \Id_{\dimp}) \) is diagonal
with the \( i \)th diagonal entry
\begin{EQA}
	\sum_{j=1}^{\dimp} \supH_{i} \supH_{j} \E (\gauss_{i} \gauss_{j} - \delta_{i,j})^{2}
	&=&
	\sum_{j=1}^{\dimp} \supH_{i} \supH_{j} + \supH_{i}^{2}
\label{8vcjv63hfgiogi3wdefjw2}
\end{EQA}
yielding
\begin{EQA}
	\E \| \BBH^{1/2} (\gaussv \gaussv^{\T} - \Id_{\dimp}) \BBH^{1/2} \|_{\Fr}^{2}
	&=&
	\sum_{i,j=1}^{\dimp} \supH_{i} \supH_{j} \E (\gauss_{i} \gauss_{j} - \delta_{i,j})^{2}
	=
	\left( \sum_{i=1}^{\dimp} \supH_{i} \right)^{2} + \sum_{i=1}^{\dimp} \supH_{i}^{2} 
\label{ikcduywjwsiv98emdvuw}
\end{EQA}
and assertion \eqref{njt66777888723fdgyf} follows.
\end{proof}

Now we compute the exponential moments of centered and non-centered quadratic forms.

\begin{lemma}
\label{Lqfexpmom}
Let \( \| \BBH \|_{\oper} = \supH \) and \( \gaussv \sim \ND(0,\Id_{\dimp}) \).
Then for any \( \mu \in (0,\supH^{-1}) \), 
\begin{EQA}
	\E \exp \Bigl\{ \frac{\mu}{2} \langle \BBH \gaussv, \gaussv \rangle \Bigr\}
	&=&
	\det(\Id_{\dimp} - \mu \BBH)^{-1/2} \, .
\label{m2v241m41m}
\end{EQA}
Moreover, with \( \dimH = \tr \BBH \) and \( \vH^{2} = \tr \BBH^{2} \)
\begin{EQA}
	\log \E \exp \Bigl\{ \frac{\mu}{2} \bigl( \langle \BBH \gaussv, \gaussv \rangle - \dimH \bigr) \Bigr\}
	& \leq &
	\frac{\mu^{2} \vH^{2}}{4 (1 - \supH \mu)} \, .
\label{m2v241m41mb}
\end{EQA}
If \( \BBH \) is positive semidefinite, \( \supH_{j} \geq 0 \), then 
\begin{EQA}
	\log \E \exp \Bigl\{ - \frac{\mu}{2} \bigl( \langle \BBH \gaussv, \gaussv \rangle - \dimH \bigr) \Bigr\}
	& \leq &
	\frac{\mu^{2} \vH^{2}}{4} \, .
\label{m2v241m41mbn}
\end{EQA}
For any complex valued \( \muH \) with \( |\supH \muH| < 1 \),
\begin{EQA}
	\biggl| \log \E \exp \Bigl\{ 
			\frac{\mu}{2} \bigl( \langle \BBH \gaussv, \gaussv \rangle - \dimH \bigr) - \frac{\muH^{2} \tr \BBH^{2}}{4}
		\Bigr\} 
	\biggr|
	& \leq &
	\frac{\supH |\mu|^{3} \vH^{2} }{6 (1 - \supH |\mu|)} \, .
\label{vbu7j3hg8hryhghyidegwdg}
\end{EQA}
\end{lemma}

\begin{proof}
W.l.o.g. assume \( \supH = 1 \).
Let \( \supH_{j} \) be the eigenvalues of \( \BBH \), 
\( |\supH_{j}| \leq 1 \).
By an orthogonal transform, one can reduce the statement to the case of a diagonal matrix 
\( \BBH = \diag\bigl( \supH_{j} \bigr) \). 
Then \( \langle \BBH \gaussv, \gaussv \rangle = \sum_{j=1}^{\dimp} \supH_{j} \gauss_{j}^{2} \) and 
by independence of the \( \gauss_{j} \)'s
\begin{EQA}
	&& \nquad
	\E \Bigl\{ \frac{\mu}{2} \langle \BBH \gaussv, \gaussv \rangle  \Bigr\}
	=
	\prod_{j=1}^{\dimp} \E \exp \Bigl( \frac{\mu}{2} \supH_{j} \eps_{j}^{2} \Bigr)
	=
	\prod_{j=1}^{\dimp} \frac{1}{\sqrt{1 - \mu \supH_{j}}} 
	=
	\det \bigl( \Id_{\dimp} - \mu \BBH \bigr)^{-1/2} .
\label{dOImuBm12EB}
\end{EQA}
Below we use the simple bound: 
\begin{EQ}[rcl]
\label{lo1uusk2iukkp}
	- \log(1 - u) - u
	&=&
	\sum_{k=2}^{\infty} \frac{u^{k}}{k}
	\leq 
	\frac{u^{2}}{2} \sum_{k=0}^{\infty} u^{k} 
	=
	\frac{u^{2}}{2 (1 - u)} \, ,
	\qquad 
	u \in (0,1),
	\\
	- \log(1 - u) + u
	&=&
	\sum_{k=2}^{\infty} \frac{u^{k}}{k}
	\leq 
	\frac{u^{2}}{2} \, ,
	\qquad \qquad
	u \in (-1,0).
\label{lo1uusk2iukk}
\end{EQ}
Now it holds for \( \mu > 0 \)
\begin{EQA}
	&& \nquad
	\log \E \Bigl\{ \frac{\mu}{2} \bigl( \langle \BBH \gaussv, \gaussv \rangle - \dimH \bigr) \Bigr\}
	=
	\log \det(\Id_{\dimp} - \mu \BBH)^{-1/2} - \frac{\mu \, \dimH}{2}
	\\
	&=&
	- \frac{1}{2} \sum_{j=1}^{\dimp} \bigl\{ \log(1 - \mu \supH_{j}) + \mu \supH_{j} \bigr\}
	\leq 
	\sum_{j=1}^{\dimp} \frac{\mu^{2} \supH_{j}^{2}}{4 (1 - \mu \supH_{j})} 
	\leq 
	\frac{\mu^{2} \vH^{2}}{4 (1 - \mu \supH)} \, .
\label{m2v241m4mj1pd}
\end{EQA}
Similarly for any complex \( \mu \) with \( |\mu| \supH < 1 \)
\begin{EQA}
	&& \nquad
	\left| 
		\log \E \Bigl\{ \frac{\mu}{2} \bigl( \langle \BBH \gaussv, \gaussv \rangle - \dimH \bigr) 
		- \frac{\muH^{2} \tr \BBH^{2}}{4} \Bigr\}
	\right|
	=
	\left| \log \det(\Id_{\dimp} - \mu \BBH)^{-1/2} - \frac{\mu \, \dimH}{2} - \frac{\muH^{2} \tr \BBH^{2}}{4} \right|
	\\
	&=&
	\frac{1}{2} \left| 
		\sum_{j=1}^{\dimp} \biggl\{ \log(1 - \mu \supH_{j}) - \mu \supH_{j} - \frac{\mu^{2} \supH_{j}^{2}}{2} \biggr\} 
	\right|
	\leq 
	\sum_{j=1}^{\dimp} \frac{|\mu \supH_{j}|^{3}}{6 (1 - |\mu|)} 
	=
	\frac{| \mu |^{3} \supH \vH^{2}}{6 (1 - |\mu|)} \, .
\label{m2v241m4mj1pd}
\end{EQA}
Statement \eqref{m2v241m41mbn} can be proved similarly.
\end{proof}

Now we consider the case of a non-centered quadratic form
\( \langle \BBH \gaussv,\gaussv \rangle/2 + \langle \Av,\gaussv \rangle \) for a fixed vector \( \Av \).

\begin{lemma}
\label{Lexpmomnoncen}
Let \( \| \BBH \| = \supH_{\max}(\BBH) < 1 \). 
Then for any \( \Av \)
\begin{EQA}
	\E \exp\Bigl\{ \frac{1}{2}\langle \BBH \gaussv,\gaussv \rangle + \langle \Av,\gaussv \rangle \Bigr\}
	&=&
	\exp\Bigl\{ \frac{\| (\Id_{\dimp} - \BBH)^{-1/2} \Av \|^{2}}{2} \Bigr\} \, \det(\Id_{\dimp} - \BBH)^{-1/2} .
\label{EeBf12BggA}
\end{EQA}
Moreover, for any \( \mu \in (0,1) \)
\begin{EQA}
	&& \nquad
	\log \E \exp\Bigl\{ 
		\frac{\mu}{2} \bigl( \langle \BBH \gaussv,\gaussv \rangle - \dimH \bigr) + \langle \Av,\gaussv \rangle 
	\Bigr\}
	\\
	&=&
	\frac{\| (\Id_{\dimp} - \mu \BBH)^{-1/2} \Av \|^{2}}{2} + \log \det(\Id_{\dimp} - \mu \BBH)^{-1/2} - \mu \, \dimH 
	\\
	& \leq &
	\frac{\| (\Id_{\dimp} - \mu \BBH)^{-1/2} \Av \|^{2}}{2} + \frac{\mu^{2} \vH^{2}}{4 (1 - \mu \| \BBH \|)} \, .
\label{EeBf12BggAmu}
\end{EQA}
\end{lemma}

\begin{proof}
Denote \( \av = (\Id_{\dimp} - \BBH)^{-1/2} \Av \). 
It holds by change of variables \( (\Id_{\dimp} - \BBH)^{1/2} \xv = \uv \) for \( \CONSTi_{\dimp} = (2\pi)^{-\dimp/2} \)
\begin{EQA}
	&& \nquad
	\E \exp\Bigl\{ \frac{1}{2}\langle \BBH \gaussv,\gaussv \rangle + \langle \Av,\gaussv \rangle \Bigr\}
	=
	\CONSTi_{\dimp}
	\int \exp\Bigl\{ - \frac{1}{2}\langle (\Id_{\dimp} - \BBH) \xv,\xv \rangle + \langle \Av,\xv \rangle \Bigr\} d\xv
	\\
	&=&
	\CONSTi_{\dimp}
	\det(\Id_{\dimp} - \BBH)^{-1/2}
	\int \exp\Bigl\{ - \frac{1}{2} \| \uv \|^{2} + \langle \av,\uv \rangle \Bigr\} d\uv
	=
	\det(\Id_{\dimp} - \BBH)^{-1/2} \, 	\ex^{\| \av \|^{2}/2}  	.
\label{EeBf12BggAp}
\end{EQA}
The last inequality \eqref{EeBf12BggAmu} follows by \eqref{m2v241m41mb}.
\end{proof}

\Section{Deviation bounds for Gaussian quadratic forms}
\label{SdevboundGauss}
The next result explains the concentration effect of \( \| \QP \xiv \|^{2} \)
for a centered Gaussian vector \( \xiv \sim \ND(0,\HVB^{2}) \) and a linear operator \( \QP \colon \R^{\dimp} \to \R^{\dimq} \),
\( \dimp,\dimq \leq \infty \).
We use a version from \cite{laurentmassart2000}.
For completeness, we present a simple proof.

\begin{theorem}
\label{TexpbLGA}
\label{Lxiv2LD}
\label{Cuvepsuv0}
Let \( \xiv \sim \ND(0,\HVB^{2}) \) be a Gaussian element in \( \R^{\dimp} \) and let
\( \QP \colon \R^{\dimp} \to \R^{\dimq} \) be such that \( \BBH = \QP \HVB^{2} \QP^{\T} \) 
is a trace operator in \( \R^{\dimq} \).
Then with \( \dimH = \tr(\BBH) \), \( \vH^{2} = \tr(\BBH^{2}) \), and 
\( \supH = \| \BBH \| \), it holds for each \( \xx \geq 0 \)
\begin{EQA}
\label{Pxiv2dimAvp12}
	\P\Bigl( \| \QP \xiv \|^{2} - \dimH > 2 \vH \, \sqrt{\xx} + 2 \supH \xx \Bigr)
	& \leq &
	\ex^{-\xx} ,
	\\
	\P\Bigl( \| \QP \xiv \|^{2} - \dimH \leq - 2 \vH \, \sqrt{\xx} \Bigr)
	& \leq &
	\ex^{-\xx} .
\label{Pxiv2dimAvp12m}
\end{EQA}
It also implies 
\begin{EQA}
	\P\bigl( \bigl| \| \QP \xiv \|^{2} - \dimH \bigr| > \zq_{2}(\BBH,\xx) \bigr)
	& \leq &
	2 \ex^{-\xx} ,
\label{PxivTBBdimA2vp}
\end{EQA}
with
\begin{EQA}
	\zq_{2}(\BBH,\xx)
	& \eqdef &
	2 \vH \, \sqrt{\xx} + 2 \supH \xx \,\, .
\label{zqdefGQF}
\end{EQA}
%
\end{theorem}

\begin{proof}
W.l.o.g. assume that \( \supH = \| \BBH \| = 1 \).
We use the identity \( \| \QP \xiv \|^{2} = \langle \BBH \gaussv, \gaussv \rangle \) with
 \( \gaussv \sim \ND(0,\Id_{\dimq}) \).
We apply Markov's inequality: with \( \mu > 0 \)
\begin{EQA}
	\P\Bigl( \langle \BBH \gaussv, \gaussv \rangle - \dimH > \zq_{2}(\BBH,\xx) \Bigr)
	& \leq &
	\E \exp \Bigl( \frac{\mu}{2} \bigl( \langle \BBH \gaussv, \gaussv \rangle - \dimH \bigr) - \frac{\mu \, \zq_{2}(\BBH,\xx)}{2} 
	\Bigr) \, .
\label{PBggiz2E2mz2}
\end{EQA}
Given \( \xx > 0 \), fix \( \mu < 1 \) by the equation
\begin{EQA}
	\frac{\mu}{1 - \mu} 
	&=&
	\frac{2 \sqrt{\xx}}{\vH} \, 
	\quad \text{ or } \quad
	\mu^{-1} 
	=
	1 + \frac{\vH}{2 \sqrt{\xx}} \, .
\label{1v2sxm12m1m}
\end{EQA}
Let \( \supH_{j} \) be the eigenvalues of \( \BBH \), 
\( |\supH_{j}| \leq 1 \).
It holds with \( \dimH = \tr \BBH \) in view of \eqref{m2v241m41mb}
\begin{EQA}
	&& \nquad
	\log \E \Bigl\{ \frac{\mu}{2} \bigl( \langle \BBH \gaussv, \gaussv \rangle - \dimH \bigr) \Bigr\}
	\leq 
	\frac{\mu^{2} \vH^{2}}{4 (1 - \mu)} \, .
\label{m2v241m4mj1p}
\end{EQA}
For \eqref{Pxiv2dimAvp12}, it remains to check that the choice \( \mu \) by \eqref{1v2sxm12m1m} yields
\begin{EQA}
	\frac{\mu^{2} \vH^{2}}{4 (1 - \mu)} - \frac{\mu \, \zq_{2}(\BBH,\xx)}{2}
	& = &
	\frac{\mu^{2} \vH^{2}}{4 (1 - \mu)} - \mu \bigl( \vH \sqrt{\xx} + \xx \bigr)
	=
	\mu \Bigl( \frac{\vH \sqrt{\xx}}{2} - \vH \sqrt{\xx} - \xx \Bigr)
	=
	- \xx .
\label{m2vA241muz2}
\end{EQA}
The bound \eqref{Pxiv2dimAvp12m} is obtained similarly from Markov's inequality 
applied to \( - \langle \BBH \gaussv, \gaussv \rangle + \dimH \) with \( \mu = 2 \vH^{-1} \sqrt{\xx} \).
The use of \eqref{m2v241m41mbn} yields
\begin{EQA}
	&& \nquad
	\P\Bigl( \langle \BBH \gaussv, \gaussv \rangle - \dimH < - 2 \vH \sqrt{\xx} \Bigr)
	\leq
	\E \exp \Bigl\{ \frac{\mu}{2} \bigl( - \langle \BBH \gaussv, \gaussv \rangle + \dimH \bigr) - \mu \, \vH \sqrt{\xx} 
	\Bigr\}
	\\
	& \leq &
	\exp \Bigl( \frac{\mu^{2} \vH^{2}}{4} - \mu \, \vH \sqrt{\xx} \Bigr) 
	=
	\ex^{-\xx} \, 
\label{PBggiz2E2mz2}
\end{EQA}
as required.
\end{proof}

\begin{corollary}
\label{CTexpbLGAd}
Assume the conditions of Theorem~\ref{TexpbLGA}.
Then for \( \zq > \vH \)
\begin{EQA}
	\P\bigl( \bigl| \| \QP \xiv \|^{2} - \dimH \bigr| \ge \zq \bigr)
	& \leq &
	2 \exp\biggl\{ - \frac{\zq^{2}}{\bigl( \vH + \sqrt{\vH^{2} + 2 \supH \zq} \bigr)^{2}} \biggr\}
	\leq 
	2 \exp\biggl( - \frac{\zq^{2}}{4\vH^{2} + 4 \supH \zq} \biggr) .
	\qquad
	\qquad
\label{3z2spsp2z3z2}
\end{EQA}
\end{corollary}

\begin{proof}
Given \( \zq \), define \( \xx \) by 
\( 2 \vH \sqrt{\xx} + 2 \supH \xx = \zq \) or 
\( 2 \supH \sqrt{\xx} = \sqrt{\vH^{2} + 2 \supH \zq} - \vH \).
Then
\begin{EQA}
	\P\bigl( \| \QP \xiv \|^{2} - \dimH \ge \zq \bigr)
	& \leq &
	\ex^{-\xx} 
	=
	\exp\biggl\{ - \frac{\bigl( \sqrt{\vH^{2} + 2 \supH \zq} - \vH \bigr)^{2}}{4 \supH^{2}} \biggr\}
	=
	\exp\biggl\{ - \frac{\zq^{2}}{\bigl( \vH + \sqrt{\vH^{2} + 2 \supH \zq} \bigr)^{2}} \biggr\}.
\label{3emzmsp22z2c}
\end{EQA}
This yields \eqref{3z2spsp2z3z2} by direct calculus.
\end{proof}

Of course, bound \eqref{3z2spsp2z3z2} is sensible only if \( \zq \gg \vH \).

\begin{corollary}
\label{RsochpHsA}
Assume the conditions of Theorem~\ref{TexpbLGA}.
If also \( \BBH \geq 0 \), then 
\begin{EQA}
\label{Pxiv2dimAxx12}
	\P\Bigl( \| \QP \xiv \|^{2} \geq \zq^{2}(\BBH,\xx) \Bigr)
	& \leq &
	\ex^{-\xx} 
\end{EQA}
with 
\begin{EQA}
	\zq^{2}(\BBH,\xx)
	& \eqdef &
	\dimH + 2 \vH \, \sqrt{\xx} + 2 \supH \xx
	\leq 
	\bigl( \sqrt{\dimH} + \sqrt{2 \supH \xx} \bigr)^{2} \, .
\label{zzxxppdBlroBB}
\end{EQA}
Also
\begin{EQA}
	\P\Bigl( \| \QP \xiv \|^{2} - \dimH < - 2 \vH \, \sqrt{\xx} \Bigr)
	& \leq &
	\ex^{-\xx} .
\label{Pxiv2dimAvp12d}
\end{EQA}
\end{corollary}

\begin{proof}
The definition implies \( \vH^{2} \leq \dimH \supH \).
One can use a sub-optimal choice of the value 
\( \mu(\xx) = \bigl\{ 1 + 2 \sqrt{\supH \dimH/\xx} \bigr\}^{-1} \) yielding the statement of the corollary.
\end{proof}

As a special case, we present a bound for the chi-squared distribution 
corresponding to \( \QP = \HVB^{2} = \Id_{\dimp} \), \( \dimp < \infty \).
Then \( \BBH = \Id_{\dimp} \), \( \tr (\BBH) = \dimp \), \( \tr(\BBH^{2}) = \dimp \) and \( \supH(\BBH) = 1 \).

\begin{corollary}
\label{Cchi2p}
Let \( \gaussv \) be a standard normal vector in \( \R^{\dimp} \).
Then for any \( \xx > 0 \)
\begin{EQA}[ccl]
\label{Pxi2pm2px}
	\P\bigl( \| \gaussv \|^{2} \geq \dimp + 2 \sqrt{\dimp \, \xx} + 2 \xx \bigr)
	& \leq &
	\ex^{-\xx},
	\\
	\P\bigl( \| \gaussv \| \,\,  \geq \sqrt{\dimp} + \sqrt{2 \xx} \bigr)
	& \leq &
	\ex^{-\xx} ,
\label{Pxi2pm2px12}
	\\
	\P\bigl( \| \gaussv \|^{2} \leq \dimp - 2 \sqrt{\dimp \, \xx} \bigr)
	& \leq &
	\ex^{-\xx}	.
\label{Pxi2pm2px22}
\end{EQA}
\end{corollary}

The bound of Theorem~\ref{TexpbLGA} 
can be represented as a usual deviation bound.

\begin{theorem}
\label{CTexpbLGA}
Assume the conditions of Theorem~\ref{TexpbLGA}.
For \( \yy > 0 \), define
\begin{EQA}
	\xx(\yy)
	& \eqdef &
	\frac{(\sqrt{\yy + \dimH} - \sqrt{\dimH})^{2}}{4 \supH} \, .
\label{iuvfiiow3kboieheuf}
\end{EQA}
Then
\begin{EQA}
	\P\bigl( \| \QP \xiv \|^{2} \ge \dimH + \yy \bigr)
	& \leq &
	\ex^{- \xx(\yy)} ,
\label{3emzmsp22z2}
	\\
	\E \bigl\{ (\| \QP \xiv \|^{2} - \dimH) \Ind\bigl( \| \QP \xiv \|^{2} \ge \dimH + \yy \bigr) \bigr\}
	& \leq &
	2 \Bigl( \frac{\yy + \dimH}{\supH \, \xx(\yy)} \Bigr)^{1/2} \, \, 
	\ex^{- \xx(\yy)} \, .
	\qquad
	\quad
\label{3emzmsp22z2e}
\end{EQA}
Moreover, let \( \muH > 0 \) fulfill \( \rexH =  \muH \supH + \muH \sqrt{\supH \dimH / \xx(\yy)} < 1 \). 
Then 
\begin{EQA}
	\E \bigl\{ \ex^{\muH (\| \QP \xiv \|^{2} - \dimH)/2} \Ind( \| \QP \xiv \|^{2} \ge \dimH + \yy) \bigr\}
	& \leq &
	\frac{1}{1 - \rexH} \, \exp\{ - (1 - \rexH) \xx(\yy) \} \, .
	\qquad
\label{llkknbononjm9hig4e}
\end{EQA}
\end{theorem}

\begin{proof}
Normalizing by \( \supH \) reduces the statements to the case with \( \supH = 1 \).
Define \( \eta = \| \QP \xiv \|^{2} - \dimH \) 
and
\begin{EQA}
	\zq(\xx)
	&=&
	2 \sqrt{\dimH \, \xx} + 2 \xx .
\label{0kmuy765433udgswhhh}
\end{EQA}
Then by \eqref{Pxiv2dimAvp12} \( \P(\eta \geq \zq(\xx)) \leq \ex^{-\xx} \).
Inverting the relation \eqref{0kmuy765433udgswhhh} yields
\begin{EQA}
	\xx(\zq)
	&=&
	\frac{1}{4} \bigl( \sqrt{\zq + \dimH} - \sqrt{\dimH} \bigr)^{2}
\label{jkv78fdjryfgsdfghgj}
\end{EQA}
and \eqref{3emzmsp22z2} follows by applying \( \zq = \yy \).
Further, 
\begin{EQA}
	\E \bigl\{ \eta \Ind(\eta \geq \yy) \bigr\}
	&=&
	\int_{\yy}^{\infty} \P(\eta \geq \zq) \, d\zq
	\leq 
	\int_{\yy}^{\infty} \ex^{ - \xx(\zq) } \, d\zq
	= 
	\int_{\xx(\yy)}^{\infty} \ex^{-\xx} \, \zq'(\xx) \, d\xx \, .
\label{zEe2Iezz2c2H23}
\end{EQA} 
As \( \zq'(\xx) = 2 + \sqrt{\dimH/\xx} \) monotonously decreases with \( \xx \), we derive
\begin{EQA}
	\E \bigl\{ \eta \Ind(\eta \geq \yy) \bigr\}
	& \leq &
	\zq'(\xx(\yy)) \ex^{-\xx(\yy)}
	=
	\frac{1}{\xx'(\yy)} \, \ex^{- \xx(\yy)}
	=
	\frac{4 \sqrt{\yy + \dimH}}{\sqrt{\yy + \dimH} - \sqrt{\dimH}} \, \ex^{- \xx(\yy)}
\label{e7ygv76bgughytuj}
\end{EQA}
and \eqref{3emzmsp22z2e} follows.

In a similar way, define \( \zqe(\xx) \) from the relation
\( \muH^{-1} \log \zqe(\xx) = \sqrt{\dimH \, \xx} + \xx \) yielding
\begin{EQA}
	\zqe(\xx)
	&=&
	\exp \bigl( \muH \sqrt{\dimH \, \xx} + \muH \, \xx \bigr) .
\label{jvcjjuvue37r6gtur4r}
\end{EQA}
The inverse relation reads
\begin{EQA}
	\xxe(\zqe)
	&=&
	\bigl( \sqrt{\muH^{-1} \log \zqe + \dimH/4} - \sqrt{\dimH/4} \bigr)^{2} .
\label{jkv78fdjryfgsdfghgjex}
\end{EQA}
Then with \( \xx(\yy) = \xxe(\ex^{\muH \yy/2}) = \bigl( \sqrt{\yy + \dimH} - \sqrt{\dimH} \bigr)^{2}/4 \)
\begin{EQA}
	\E \bigl\{ \ex^{\muH \eta/2} \Ind(\eta \geq \yy) \bigr\}
	&=&
	\int_{\ex^{\muH \yy/2}}^{\infty} \P(\ex^{\muH \eta/2} \geq \zqe) \, d\zqe
	=
	\int_{\ex^{\muH \yy/2}}^{\infty} \P(\eta \geq 2\muH^{-1} \log \zqe) \, d\zqe
	\\
	& \leq &
	\int_{\ex^{\muH \yy/2}}^{\infty} \ex^{ - \xxe(\zqe) } \, d\zqe
	= 
	\int_{\xx(\yy)}^{\infty} \ex^{-\xx} \, \zqe'(\xx) \, d\xx .
\label{zEe2Iezz2c2H23}
\end{EQA} 
Further, in view of \( \muH + 0.5 \,\muH \sqrt{\dimH/\xx} < \muH + \muH \sqrt{\dimH / \xx(\yy)} = \rexH < 1 \) for 
\( \xx \geq \xx(\yy) \), it holds
\begin{EQA}
	\zqe'(\xx)
	&=&
	\bigl( \muH + 0.5 \, \muH \sqrt{\dimH/\xx} \bigr) \exp \bigl( \muH \sqrt{\dimH \, \xx} + \muH \, \xx \bigr) 
	\leq 
	\exp \bigl( \muH \, \xx \sqrt{\dimH / \xx(\yy)} + \muH \, \xx \bigr)
	=
	\exp (\rexH \, \xx) 
\label{jcuyu3ww3jbkihjitwedk}
\end{EQA}
and  
\begin{EQA}
	\E \bigl\{ \ex^{\muH \eta/2} \Ind(\eta \geq \yy) \bigr\}
	& \leq &
	\int_{\xx(\yy)}^{\infty} \ex^{-(1 - \rexH)\xx} \, d\xx 
	=
	\frac{1}{1 - \rexH} \, \ex^{-(1 - \rexH)\xx(\yy)} \, 
\label{zEe2Iezz2c2H23}
\end{EQA} 
and \eqref{llkknbononjm9hig4e} follows.
\end{proof}

\bibliography{exp_ts,listpubm-with-url}

\end{document}